\documentclass[12pt,draft,leqno]{article}
\usepackage{amssymb,eucal,latexsym}
\usepackage{amsmath,amsthm,latexsym,amscd,amsbsy,fleqn,graphicx}

\newtheorem{theorem}{Theorem}[section]
\newtheorem{lm}[theorem]{Lemma}
\newtheorem{exa}[theorem]{Example}

\newtheorem{cor}[theorem]{Corollary}
\newtheorem{pro}[theorem]{Proposition}
\newtheorem{defi}[theorem]{Definition}

\newtheorem{nota}[theorem]{Notation}

\newtheorem{rem}[theorem]{Remark}

\newtheorem{nist}[theorem]{}

\def\p{\varphi}
\def\a{\alpha}

\def\d{\delta}

\def\g{\gamma}
\def\GA{\Gamma}

\def\k{\kappa}

\def\l{\lambda}

\def\s{\sigma}

\def\fs{\hat{f}}
\def\ps{\hat{\varphi}}

\def\lra{\longrightarrow}

\def\sbe{\subseteq}
\def\spe{\supseteq}
\def\stm{\setminus}
\def\ems{\emptyset}
\def\nes{\neq\emptyset}

\def\unl{\underline}
\def\unlx{\underline{X}}
\def\unlxc{\underline{X}^c}
\def\unlb{\underline{B}}
\def\unlbc{\underline{B}^c}

\def\ex{\exists}
\def\fa{\forall}

\def\we{\wedge}

\def\ap{^{\prime}}
\def\inv{^{-1}}
\def\st{\ |\ }

\def\nin{\not\in}

\def\card #1{\vert #1 \vert}

\def\AA{{\cal A}}
\def\BB{{\cal B}}
\def\CC{{\cal C}}

\def\FF{{\cal F}}
\def\GG{{\cal G}}

\def\PP{{\cal P}}

\def\RR{{\cal R}}
\def\SS{{\cal S}}
\def\TT{{\cal T}}

\def\PCR{{\rm PCRel}}
\def\CR{{\rm CRel}}
\def\CLR{{\rm CloRel}}
\def\CLRSR{{\rm CloRSRel}}

\def\Stone{{\bf Stone}}

\def\CBool{{\bf CBool}}

\def\1{{\bf 1}}
\def\2{\mbox{{\bf 2}}}
\def\3{\mbox{{\bf 3}}}

\def\int{\mbox{{\rm int}}}

\def\cl{\mbox{{\rm cl}}}

\def\doc{\hspace{-1cm}{\em Proof.}~~}
\def\sq{\hspace*{\fill} \hbox{\vrule\vbox{\hrule\phantom{o}\hrule}\vrule}}
\def\sqs{\sq \vspace{2mm}}

\def\Boo{{\bf Bool}}
\def\ECS{{\bf 2Stone}}

\def\ECC{{\bf ECS}}
\def\CStone{{\bf 2CStone}}

\def\tcx{t_X^C}
\def\tcy{t_Y^C}
\def\tcx0{t_{(X,X_0)}}
\def\tcy0{t_{(Y,Y_0)}}

\def\CSR{\textsl{CSR}}
\def\CCSR{\textsl{CCSR}}

\def\CRel{\textsl{CRel}}
\def\CCRel{\textsl{CCRel}}

\def\wh{\widehat}

\def\bU0{\bar{U}=(U^0,(U^i,U^{ci})_{i\in\omega})}

\def\bV0{\bar{V}=(V^0,(V^i,V^{ci})_{i\in\omega})}

\def\PCA{{\cal PCA}}
\def\PCS{{\cal PCS}}

\title{{\LARGE\bf Topological Representation}\\
\vspace{0.2cm}
{\LARGE\bf of Precontact Algebras}\\
\vspace{0.2cm}
{\LARGE\bf and a Connected Version of}\\
\vspace{0.2cm}
{\LARGE\bf  the Stone Duality Theorem -- I}\\
\vspace{0.5cm}
{\large\bf Georgi Dimov\thanks{{\footnotesize The first (resp., the second) author of this paper was supported by the
contract no. 7/2015 $``$Contact algebras and extensions of topological spaces" (resp., contract no. 5/2015) of the
Sofia University Science Fund.}}\hspace{1mm}
 and Dimiter Vakarelov}\\
\vspace{0.2cm}
{\footnotesize\rm Department of Mathematics and Informatics,  University of Sofia,}\\
{\footnotesize\rm 5 J. Bourchier Blvd., 1164 Sofia, Bulgaria}
}

\author{}

\date{}

\begin{document}

\maketitle

\begin{abstract}
The notions of a {\em 2-precontact space}\/ and a {\em 2-contact space}\/
are introduced. Using
them, new representation theorems for  precontact and contact
algebras
are proved. They
incorporate and strengthen both the discrete and topological
representation theorems from \cite{DV,DV1,DV2,DW,VDDB}. It is
shown that there are bijective correspondences between such kinds
of algebras and such kinds of spaces.
As applications of the obtained results, we get new  connected versions of the Stone Duality Theorems \cite{ST,Si}  for Boolean algebras and for complete Boolean algebras,
as well as a Smirnov-type theorem (in the sense of \cite{Sm})
for a kind of compact $T_0$-extensions of compact Hausdorff extremally disconnected spaces.
We also introduce the notion of a
{\em Stone
adjacency space}\/  and using it, we prove another representation
theorem for precontact algebras.  We even obtain a  bijective correspondence between the class of
all, up to isomorphism, precontact algebras and the class of
all, up to iso\-mor\-ph\-ism, Stone adjacency spaces.
\end{abstract}



\footnotetext[3]{{\footnotesize {\em E-mail addresses:}
gdimov@fmi.uni-sofia.bg, dvak@fmi.uni-sofia.bg}}

\section{Introduction}
In this paper we give the proofs of the results announced in the first seven sections of our paper \cite{DV3} (the results of the remaining sections will be proved in the second part of the paper) and obtain many new additional results and some new applications.
In it we present a common approach both to the discrete
and to the non-discrete region-based theory of space. The paper is a
continuation of the line of investigations started in \cite{VDDB}
and continued in \cite{DV1,DV2,DW}.

 Standard models of non-discrete theories of space are the
contact algebras of regular closed subsets of some topological
spaces (\cite{Stell,Ro,VDDB,DV1,DV2,DW}). In a sense these topological
models reflect the continuous nature of the space. However, in the
``real-world" applications, where digital methods of modeling are
used, the continuous models of space are not enough.
This motivates   a search for good ``discrete" versions of the
theory of space. One kind of discrete models are the so called
{\em adjacency spaces}, introduced by Galton \cite{Galton} and
generalized by D\"{u}ntsch and Vakarelov in \cite{DV}. Based on
the Galton's approach, Li and Ying \cite{Li} presented a
``discrete" generalization of the Region Connection Calculus
(RCC). The latter, introduced in \cite{R}, is one of the main
systems in the non-discrete region-based theory of space. A
natural class of Boolean algebras related to adjacency spaces are
the {\em precontact algebras}, introduced in \cite{DV} under the
name of {\em proximity algebras}. The notion of precontact algebra
is a generalization of the notion of contact algebra. Each
adjacency space generates canonically a precontact algebra. It is
proved in \cite{DV} (using another terminology) that each
precontact algebra
 can be embedded in the precontact algebra of an
 adjacency space.
 In \cite{DV1} we proved that each contact algebra can be
  embedded in the standard contact algebra of a compact semiregular
  $T_{0}$-space, answering the question of D\"{u}ntsch and
  Winter,
  posed in \cite{DW}, whether the contact algebras have a topological representation.
    This shows that contact algebras
  possess both a discrete  and a
  non-discrete (topological) representation. In this paper we extend the representation
  techniques developed in \cite{DV1,DV2} to precontact algebras,
  proving that each precontact algebra can be
  embedded in a special topological object, called a {\em
  2-precontact space}. We also establish a bijective correspondence
  between  precontact algebras and
  2-precontact spaces. This result is new even in the special
  case of contact algebras: introducing the notion of {\em
  2-contact space}\/ as a specialization of the notion of a 2-precontact
  space, we show that there is a bijective correspondence
  between  contact algebras and
  2-contact spaces. Also, we introduce the notion of a
  {\em Stone
adjacency space}\/  and using it, we prove another representation
theorem for precontact algebras. We even obtain a  bijective correspondence between the class of
all, up to isomorphism, precontact algebras and the class of
all, up to iso\-mor\-ph\-ism, Stone adjacency spaces.

The developed theory permits us to obtain as  corollaries the celebrated Stone Representation Theorem \cite{ST}
and a new connected version of it.  They correspond, respectively, to the extremal contact relations on Boolean algebras: the smallest one and the largest one.
We show as well that the new connected version of the Stone Representation Theorem can be extended to a new connected version of the Stone Duality Theorem.
Let us explain what we mean by a $``$connected version".
The celebrated Stone Duality Theorem \cite{ST} states that the category $\Boo$ of all Boolean algebras and Boolean homomorphisms is dually equivalent to the category $\Stone$ of  compact Hausdorff totally disconnected spaces (i.e., {\em Stone spaces}\/) and continuous maps.  The restriction of the Stone duality to the category $\CBool$ of  complete Boolean algebras and Boolean homomorphisms is a duality between the category $\CBool$ and the category of   compact Hausdorff extremally disconnected spaces  and continuous maps. We introduce the notion of a {\em Stone 2-space} and the category $\ECS$ of Stone 2-spaces and suitable morphisms between them, and we show that the category $\ECS$ is dually equivalent to the category $\Boo$. The Stone 2-spaces are pairs $(X,X_0)$ of a compact {\em connected} $T_0$-space $X$ and a dense subspace $X_0$ of $X$, satisfying some mild conditions. We introduce as well the notion of an {\em extremally connected space} and show that the category $\ECC$ of  extremally connected spaces and  continuous maps between them  satisfying a natural condition, is dually equivalent to the category $\CBool$. The extremally connected spaces are compact {\em connected} $T_0$-spaces satisfying an additional condition, and the open maps are part of the  morphisms of the category $\ECC$.

As another application of the obtained results, we prove a Smirnov-type theorem (in the sense of \cite{Sm}). In his celebrated Compactification Theorem, Ju. M. Smirnov \cite{Sm} proved that there
exists an isomorphism between the ordered set of all Efremovi\v{c} proximities on a
Tychonoff space X and all, up to equivalence, compact Hausdorff extensions of X.
The notion of a contact relation on a Boolean algebra is a generalization of the notion
of a proximity. We   show
that there exist an isomorphism between the ordered by  inclusion set of all contact
relations on a complete Boolean algebra $B$ and the ordered (by the injective order) set of all, up to isomorphism, {\em C-semiregular}\/ extensions
of its Stone space $S(B)$. In this way we describe C-semiregular extensions of extremally disconnected compact Hausdorff spaces. The notion of a {\em C-semiregular space}\/ was introduced in \cite{DV1}. It appears naturally in the theory of contact algebras. The class of C-semiregular spaces is a subclass of the class of compact $T_0$-spaces. As a corollary, we obtain that every extremally disconnected compact Hausdorff space $X$ has a largest C-semiregular extension $(\g X,\g_X)$; moreover, $\g X$ is an extremally connected space and this characterizes it between all C-semiregular extensions of $X$. We show that every continuous map $f:X\lra Y$ between two extremally disconnected compact Hausdorff spaces $X$ and $Y$ has a continuous extension $\g f:\g X\lra \g Y$. We obtain, as well, some other similar results about continuous extensions of continuous maps.


The paper is organized as follows. In Section 2 we
introduce the notions of  precontact and contact algebra and give
the two main examples of them: the precontact algebras on
adjacency spaces, and the contact algebras on topological spaces.
In Section 3 we introduce three kinds of points in precontact
algebras: ultrafilters, grills and clans.
Also, the notions of a {\em topological adjacency space}\/ and a {\em Stone
adjacency space}\/ are introduced and our first representation
theorem for precontact algebras is proved there. In Section 4 we
introduce the notions of  {\em 2-precontact space} and  {\em
canonical precontact algebra of a 2-precontact space}. In Section
5 we associate with each precontact algebra $\unlb$ a
2-precontact space, called the {\em canonical 2-precontact space
of}\/ $\unlb$. In Section 6 we present the main theorem of
the paper: the second representation theorem for precontact algebras. In
Section 7 we introduce the notion of a {\em 2-contact space} and we
prove that there exists a bijective correspondence between the class of all (up
to isomorphism) contact algebras and the class of all, up to isomorphism,
2-contact spaces. This is a generalization of the similar result
about complete contact algebras obtained in \cite{DV1}. In Section
8, the results of which are completely new and were not announced in our paper \cite{DV3}, we
demonstrate that the  Stone Representation Theorem \cite{ST}
and the new connected version of it, which we now obtain, follow from our representation theorem for contact algebras presented in Section 7.
Here we obtain also the new connected versions of the Stone Duality Theorems for Boolean algebras and for complete Boolean algebras, about which we already mentioned above.
In the last Section 9, we collect our results about C-semiregular extensions of extremally disconnected compact Hausdorff spaces about which we also mentioned above.
These results  were  not presented in the paper \cite{DV3} and are new.

 We now fix the notations.

 All lattices are with top (= unit) and bottom (= zero) elements,
denoted respectively by 1 and 0. We do not require the elements
$0$ and $1$ to be distinct.

 If $(X,\tau)$ is a topological space and $M$ is a subset of $X$, we
denote by $\cl_{(X,\tau)}(M)$ (or simply by $\cl(M)$ or
$\cl_X(M)$) the closure of $M$ in $(X,\tau)$ and by
$\int_{(X,\tau)}(M)$ (or briefly by $\int(M)$ or $\int_X(M)$) the
interior of $M$ in $(X,\tau)$. The open maps between topological spaces are supposed to be continuous.
The extremally disconnected spaces and compact spaces are not assumed to be Hausdorff (as it is adopted in \cite{E}).

If $X$ is a topological space, we denote by $CO(X)$ the set of
all clopen subsets of $X$. Obviously, $(CO(X),\cup,\cap,\stm,\ems, X)$ is a Boolean algebra.

 If $X$ is a set, we denote by $2^X$ the power set of $X$.

 If $\CC$ denotes a category, we write $X\in \card\CC$ if $X$ is an
object of $\CC$, and $f\in \CC(X,Y)$ if $f$ is a morphism of $\CC$
with domain $X$ and codomain $Y$.

   The main reference books for all notions which are not defined here
 are \cite{E,kop89,AHS}.

\section{ Precontact algebras}

\subsection*{Precontact and contact algebras.}

\begin{defi}\label{precontact}
\rm
  An algebraic system\/ $\unlb=(B,C)$ is called a\/ {\em precontact
algebra} (\cite{DV}) (abbreviated as PCA) if the following holds:
\begin{description}
\item[$\bullet$] \quad $B=(B,0,1,+,.,*)$ is a Boolean algebra
(where the complement is denoted by $``*$");

\item[$\bullet$]\quad  $C$ is a binary relation on $B$ (called
a {\em precontact relation}\/) satisfying the following
axioms:

\item[$(C0)$]\quad \ If $aCb$ then $a\not=0$ and $b\not=0$;

\item[$(C+)$]\quad $aC(b+c)$ iff $aCb$ or $aCc$;  $(a+b)Cc$ iff
$aCc$ or $bCc$.
\end{description}

\noindent A precontact algebra $(B,C)$ is said to be\/ {\em
complete} if the Boolean algebra $B$ is complete.
 Two precontact algebras\/ $\unlb=(B,C)$ and\/
$\underline{B_1}=(B_1,C_1)$ are said to be\/ {\em PCA-isomorphic}
(or, simply,\/ {\em isomorphic}) if there exists a {\em PCA-isomorphism} between them, i.e., a Boolean
isomorphism $\varphi: B\longrightarrow B_1$ such that, for every
$a,b\in B$, $aCb$ iff $\varphi(a)C_1 \varphi(b)$.

 The negation of the relation $C$ is denoted by $(-C)$.

 For any PCA $(B,C)$, we define a binary relation  $``\ll_C $"  on
$B$ (called {\em non-tangential inclusion})  by $$\ a \ll_C b
\leftrightarrow a(-C)b^*.$$ Sometimes we will write simply
$``\ll$" instead of $``\ll_C$".

 We will also consider precontact algebras satisfying some
additional axioms:

\item[$(C ref)$]\quad \  If $a\not=0$ then $aCa$ (reflexivity
axiom);

\item[$(C sym)$]\quad If $aCb$ then $bCa$  (symmetry axiom);

\item[$(C tr)$]\quad\quad If $a\ll_C c$ then $(\exists b)(a\ll_C
b\ll_C c)$ (transitivity axiom);

\item[$(C con)$]\quad \  If $a\not=0,1$ then $aCa^{*}$ or $a^* Ca$
(connectedness axiom).

\smallskip

 A precontact algebra $(B,C)$ is called a\/ {\em contact algebra}
 (\cite{DV1})
(and $C$ is called a\/ {\em contact relation}) if it satisfies the
axioms $(C ref)$ and $(C sym)$.
We say that two contact algebras are {\em CA-isomorphic}\/ if they are PCA-isomorphic; also, a PCA-isomorphism between two contact algebras will be called a {\em CA-isomorphism}.

 A precontact algebra $(B,C)$ is
called\/ {\em connected}  if it satisfies the axiom $(Ccon)$.
\end{defi}

The following lemma says that in every precontact algebra we can
define a contact relation.

\begin{lm}\label{lemma1} Let $(B,C)$ be a precontact algebra. Define
$$aC^{\#}b\iff ((aCb)\vee  (bCa) \vee (a.b\not=0)).$$ Then $C^{\#}$ is a
contact relation on $B$ and hence $(B,C^{\#})$ is a contact
algebra.
\end{lm}

\doc If $a\neq 0$ then $a.a=a\neq 0$ and thus $aC^{\#}a$. So, $C^{\#}$ satisfies the axiom $(C ref)$.
Further, let $aC^{\#}b$. Then there are three possibilities: (1) if $aCb$ then $bC^{\#}a$; (2) if $bCa$ then $bC^{\#}a$; (3) if $a.b\neq 0$ then $b.a\neq 0$ and thus
$bC^{\#}a$. Therefore, $C^{\#}$ satisfies the axiom $(C sym)$.
\sqs

\begin{rem}\label{ctrdiez}
\rm
We will also consider precontact algebras satisfying the following
variant of the transitivity axiom (Ctr):

\smallskip

\noindent($C tr\#$)\quad\quad If $a\ll_{C^\#} c$ then $(\exists
b)(a\ll_{C^\#} b\ll_{C^\#} c)$.

\smallskip

The axiom ($C tr\#$) is known  as the $``$Interpolation axiom".

A contact algebra $(B,C)$ is called a  {\it  normal contact algebra} (\cite{deV,F}) if it satisfies the  axiom ($Ctr\#$) and the following one:

\smallskip

\noindent(C6) If $a\not= 1$ then there exists $b\not= 0$ such that
$b(-C)a$.

\smallskip

\noindent  The notion of a
normal contact algebra was introduced by Fedorchuk \cite{F} (under
the name of $``${\em Boolean $\d$-algebra}") as an equivalent expression
of the notion of a {\em compingent Boolean algebra}\/ of de Vries \cite{deV} (see its definition below). We call
such algebras $``$normal contact algebras" because they form a
subclass of the class of contact algebras and naturally arise in
normal Hausdorff spaces.

The relations $C$ and $\ll$ are inter-definable. For example,
normal contact algebras could be equivalently defined (and exactly
in this way they were introduced (under the name of {\em
compingent Boolean algebras}) by de Vries in \cite{deV}) as a pair
of a Boolean algebra $B=(B,0,1,+,.,{}^*)$ and a binary
relation $\ll$ on $B$ subject to the following axioms:

\smallskip

\noindent ($\ll$1) $a\ll b$ implies $a\leq b$;\\
($\ll$2) $0\ll 0$;\\
($\ll$3) $a\leq b\ll c\leq t$ implies $a\ll t$;\\
($\ll$4) ($a\ll b$ and $a\ll c$) implies $a\ll b.c$;\\
($\ll$5) If  $a\ll c$ then $a\ll b\ll c$  for some $b\in B$;\\
($\ll$6) If $a\neq 0$ then there exists $b\neq 0$ such that $b\ll
a$;\\
($\ll$7) $a\ll b$ implies $b^*\ll a^*$.

\smallskip

Note that if $0\neq 1$ then the axiom ($\ll$2) follows from the
axioms ($\ll$3), ($\ll$4), ($\ll$6) and ($\ll$7).

\smallskip

Obviously, contact algebras could be equivalently defined as a
pair of a Boolean algebra $B$ and a binary relation $\ll$ on $B$
subject to the  axioms ($\ll$1)-($\ll$4) and ($\ll$7); then, clearly,  the relation $\ll$
satisfies also the axioms

\smallskip

\noindent($\ll$2') $1\ll 1$;\\
($\ll$4') ($a\ll c$ and $b\ll c$) implies $(a+b)\ll c$.

 It is not difficult to see that precontact algebras could be equivalently defined as a
pair of a Boolean algebra $B$ and a binary relation $\ll$ on $B$
subject to the  axioms ($\ll$2), ($\ll$2'), ($\ll$3), ($\ll$4) and ($\ll$4').

\smallskip

It is easy to see that axiom (C6) can be stated
equivalently in the form of  ($\ll$6).
\end{rem}

\subsection*{Examples of precontact and contact algebras}

\subsubsection*{1. Extremal contact relations}

\begin{exa}\label{extrcr}
\rm Let $B$ be a Boolean algebra. Then there exist a largest and a
smallest contact relations on $B$; the largest one, $\rho_l$ (sometimes we will write $\rho_l^B$), is
defined by $$a\rho_l b \iff (a\neq 0\mbox{ and }b\neq 0),$$ and the
smallest one, $\rho_s$ (sometimes we will write $\rho_s^B$), by $$a\rho_s b \iff a\wedge b\neq 0.$$

Note that, for $a,b\in B$, $$a\ll_{\rho_s} b \iff a\le b;$$ hence
$a\ll_{\rho_s} a$, for any $a\in B$. Thus $(B,\rho_s)$ is a normal
contact algebra.
\end{exa}

\subsubsection*{2. Precontact algebras on adjacency spaces. {\rm (Galton \cite{Galton}, D{\"u}ntsch and Vakare\-lov \cite{DV})}}
By an\/ {\em adjacency space} we mean a relational  system $(W,R)$,
where
  $W$ is a non-empty set whose elements are
called\/ {\em cells}, and $R$ is a binary relation on $W$ called
the\/ {\em adjacency relation};  the subsets of $W$ are called\/
{\em regions}.

The reflexive and symmetric closure $R^{\flat}$ of $R$ is defined
as follows:
\begin{equation}\label{rflat}
xR^{\flat}y \iff ((xRy) \vee (yRx) \vee (x=y)).
\end{equation}

  A precontact relation $C_{R}$
between the  regions of an adjacency space $(W,R)$ is defined as
follows: for every $a,b\subseteq W$,
\begin{equation}\label{CR}
aC_{R} b \mbox{ iff }(\exists x\in a)(\exists y\in b)(xRy).
\end{equation}

\begin{pro}\label{prop1}{\rm (\cite{DV})} Let $(W,R)$ be an adjacency space
and let  \  $2^W$ be the Boolean algebra of all subsets of $W$.
Then:
\begin{description}
\item[{\rm (a)}]\quad  $(2^W, C_{R})$  is a precontact algebra;

\item[{\rm (b)}]\quad  $(2^W, C_{R})$  is a contact algebra iff $R$ is a
reflexive and symmetric relation on $W$. If $R$ is a reflexive and
symmetric relation on $W$ then $C_{R}$ coincides with
$(C_{R})^{\#}$ and $C_{R^{\flat}}$;

\item[{\rm (c)}]\quad   $C_{R}$ satisfies the axiom $(C tr)$ iff $R$ is
a transitive relation  on $W$;

\item[{\rm (d)}]\quad  $C_{R}$ satisfies the axiom $(Ccon)$ iff $R$ is a
connected relation on $W$ (which means that if $x,y\in W$ and
$x\not=y$ then there is an $R$-path from $x$ to $y$ or from $y$ to
$x$).
\end{description}
\end{pro}

\begin{theorem}\label{theorem1} {\rm  (\cite{DV})}  Each
precontact algebra $(B,C)$ can be isomorphically embedded into the
precontact algebra $(2^W,C_R)$ of some adjacency space $(W,R)$.
Moreover, if  $(B,C)$ satisfies  the axiom (Cref) (respectively, (Csym); (Ctr)) then the relation $R$ is reflexive (respectively,
symmetric; transitive).
\end{theorem}

\subsubsection*{3. Contact algebras on topological spaces.}

\begin{nist}\label{rcdef}
\rm
Let $X$ be a topological space and let $RC(X)$ be the set of all
regular closed subsets of $X$ (recall that a subset $F$ of $X$ is
said to be {\em regular closed}\/ if $F=\cl(\int(F))$). Let us equip
$RC(X)$ with the following Boolean operations and  contact
relation $C_{X}$:
\begin{description}
\item[$\bullet$] \quad $a+b=a\cup b$;

\item[$\bullet$] \quad $a^{*}=\cl(X\setminus a)$;

\item[$\bullet$] \quad  $a.b=\cl(\int(a \cap b)) (= (a^{*}\cup b^{*})^{*})$;

\item[$\bullet$] \quad $0=\emptyset$, $1=X$;

\item[$\bullet$] \quad  $aC_{X}b$ iff $a\cap b\not=\emptyset$.
\end{description}
\end{nist}

The following lemma is a well-known fact.

\begin{lm}\label{rcl} Let $X$ be a topological space. Then $$(RC(X), C_{X})=(RC(X),0,1,+,.,*,C_{X})$$ is a contact
algebra.
\end{lm}

The contact algebras of the type $(RC(X), C_{X})$, where $X$ is a topological space, are called {\em standard  contact algebras}.

Recall that a space $X$ is said to be\/ {\em semiregular} if
$RC(X)$ is a closed base for $X$.
Recall as well the following definition:
if $(A,\le)$ is a poset and $B\sbe A$ then $B$ is said
to be a {\em dense subset of} $A$ if for any $a\in A\stm\{0\}$
there exists $b\in B\stm\{0\}$ such that $b\le a$; when
$(B,\le_1)$ is a poset and $f:A\lra B$ is a map, then we will say
that $f$ is a {\em dense map}\/ if $f(A)$ is a dense subset of
$(B,\le_1)$.

The following theorem answers the question, posed by D{\"u}ntsch
and Winter in \cite{DW}, whether contact algebras have a
topological representation:

\begin{theorem}\label{theorem2}{\rm (\cite{DV1})}
For each  contact algebra  $\unl{B} = (B,C)$ there exists a dense embedding $g_B$
of $B$ into a standard  contact algebra  $(RC(X,\tau),C_X)$, where $(X,\tau)$ is a  compact semiregular $T_0$-space. The
algebra $B$ is connected iff the space $X$ is connected. When $B$ is complete then the embedding $g_B$ becomes an
isomorphism between contact algebras $(B,C)$ and $(RC(X),C_X)$.
\end{theorem}

 The aim of this work is to generalize  Theorem \ref{theorem1} and Theorem
 \ref{theorem2}
in several ways: to find a topological representation of
precontact
 algebras which incorporates both the
 ``discrete"  and the ``continuous" nature of the space;
 to find
representation
 theorems  in the  style of the  Stone representation  of Boolean
 algebras instead of  embedding theorems;
to establish, again as in the Stone theory, a bijective
 correspondence between precontact algebras and the corresponding
 topological objects;
to find some new applications of the obtained results.

\section{Points in precontact algebras}

 In this section we   introduce
 three kinds of abstract points in precontact
algebras:  ultrafilters, grills and  clans.
This is done by analogy with the case of contact algebras
(see, e.g., \cite{DV1,VDDB}). We assume that the notions of a
filter and ultrafilter of a Boolean algebra are familiar.  Clans
were introduced by Thron \cite{Thron} in proximity theory. Our
definition is a lattice-theoretic  generalization of Thron's definition.

{\em The set of all ultrafilters of a Boolean algebra $B$ is
denoted by} $Ult(B)$.

\begin{defi}\label{clandef}
\rm
Let\/ $\unlb=(B,C)$ be a precontact algebra.
A non-empty subset $\Gamma$ of $B$ is called a\/ {\em  clan} if it
satisfies the following conditions:
\begin{description}
\item[ $(Clan 1)$]\quad  $0\not\in \Gamma$;

\item[ $(Clan 2)$]\quad   If $a\in\Gamma$ and $a\leq b$ then
$b\in\Gamma$;

\item[ $(Clan 3)$]\quad  If $a+b\in \Gamma$ then  $a\in \Gamma$ or
$b\in \Gamma$;

\item[ $(Clan 4)$]\quad If $a,b\in \Gamma$ then $aC^{\#}b$.
\end{description}


The set of all clans
of a precontact
algebra\/ $\unlb$ is denoted by $Clans(\unlb)$.
\end{defi}

The following lemma is  obvious:

\begin{lm}\label{clan1}  Let\/ $\unlb=(B,C)$ be a precontact algebra. Each ultrafilter of $B$ is a clan of $\unlb$ and hence
$Ult(B)\subseteq Clans(\unlb)$.
\end{lm}

Now, for any  precontact algebra $\unlb=(B,C)$, we  define  a binary relation $R_{\unlb}$  between the ultrafilters of $B$
 making the set $Ult(B)$ an   adjacency space.

\begin{defi}\label{ub}
\rm
Let\/ $\unlb=(B,C)$ be a precontact algebra and let
$U_{1},U_{2}$ be ultrafilters of $B$. We set
\begin{equation}\label{ultra}
U_{1} R_{\unlb} U_{2} \mbox{ iff }(\forall a\in U_{1})(\forall b\in
U_{2})(aCb)\ \  (\mbox{i.e., iff } U_{1}\times U_{2}\subseteq C).
\end{equation}

 The relational system $(Ult(B), R_{\unlb})$ is
called the\/ {\em canonical adjacency space of}\/  $\unlb$.

We say that $U_{1},U_{2}$ are\/ {\em connected} iff
$U_{1}(R_{\unlb})^{\flat}  U_{2}$ (see (\ref{rflat}) for the notation $R^{\flat}$).
\end{defi}

The next lemma is obvious.

\begin{lm}\label{clan2} Let\/ $\unlb=(B,C)$ be a precontact algebra and let $I$ be a  set of
 connected ultrafilters. Then the union
$\GA = \bigcup\{U\ |\ U\in I\}$ is a clan.
\end{lm}

\begin{lm}\label{lmpc}{\rm (\cite{DV1,DV})}
\textsc{(Ultrafilter and clan characterizations of precontact and
contact relations.)}
Let\/ $\unlb=(B,C)$ be a precontact algebra and $(Ult(B),R_{\unlb})$
be the canonical adjacency space of\/  $\unlb$. Then the
following is true for any $a,b\in B$:
\begin{description}
\item[{\rm (a)}]\quad $aCb$ iff $(\exists U_{1},U_{2}\in
Ult(B))((a\in U_1)\we (b\in U_2)\we (U_{1}R_{\unlb} U_{2}))$;

\item[{\rm(b)}]\quad $aC^{\#}b$ iff $(\exists U_{1},U_{2}\in
Ult(B))((a\in U_1)\we (b\in U_2)\we (U_{1}R^{\flat}_{(B,C)} U_{2}))$;

\item[{\rm(c)}]\quad $aC^{\#}b$ iff $(\exists \Gamma\in Clans(\unlb))(a,b\in \Gamma)$;

\item[{\rm(d)}]\quad $R_{\unlb}$ is a reflexive relation iff\/ $\unlb$
satisfies the axiom $(Cref)$;

\item[{\rm(e)}]\quad $R_{\unlb}$ is a symmetric relation iff\/ $\unlb$
satisfies the axiom $(Csym)$;

\item[{\rm(f)}]\quad $R_{\unlb}$ is a transitive relation iff\/ $\unlb$
satisfies the axiom $(Ctr)$.
\end{description}
\end{lm}

Recall that a non-empty subset of a Boolean algebra $B$ is called a {\em grill}\/ if it satisfies the axioms (Clan1)-(Clan3). The set of all grills of $B$ will be denoted by
$Grills(B)$.
The next lemma is well known (see, e.g., \cite{Thron}):

\begin{lm}\label{grilllemma} \textsc{(Grill Lemma.)} If $F$ is a filter of a Boolean algebra $B$ and $G$ is a grill of $B$ such that $F\sbe G$ then there exists an
ultrafilter $U$ of $B$ with $F \sbe U\sbe G$.
\end{lm}

\subsection*{Stone adjacency spaces and representation of precontact algebras}

\begin{defi}\label{f}
\rm
 Let $X$ be a non-empty topological space and $R$ be a binary relation on $X$.
  Then the pair $(CO(X), C_{R})$ (see (\ref{CR})
for $C_R$) is a precontact algebra (by Proposition \ref{prop1}(a)), called the {\em canonical precontact
algebra of the relational system $(X,R)$.}
\end{defi}

\begin{defi}\label{tas}
\rm
An adjacency space $(X,R)$ is called a\/ {\em topological
adjacency space} (abbreviated as TAS) if $X$ is a topological
space and $R$ is a closed relation on $X$. When $X$ is a compact
Hausdorff zero-dimensional space (i.e., when $X$ is a\/ {\em Stone
space}), we say that the topological adjacency space $(X,R)$ is
a\/ {\em Stone adjacency space}.

Two topological adjacency spaces $(X,R)$ and $(X_1,R_1)$ are said
to be\/ {\em TAS-isomorphic} if there exists a homeomorphism
$f:X\longrightarrow X_1$ such that, for every $x,y\in X$, $xRy$
iff $f(x)R_1 f(y)$.
\end{defi}

Recall that the Stone space
$S(A)$ of a Boolean algebra $A$ is the set $X=Ult(A)$ endowed
with a topology $\TT$ having as a closed base the family
$\{s_A(a)\st a\in A\}$, where
\begin{equation}\label{sofa}
s_A(a)=\{u\in X\st a\in u\},
\end{equation}
for every $a\in A$; then $$S(A)=(X,\TT)$$ is a compact Hausdorff
zero-dimensional space, $s_A(A)= CO(X)$ and {\em the Stone map}
\begin{equation}\label{stonemap}
s_A:A\lra CO(X), \ \ a\mapsto s_A(a),
\end{equation}
is a Boolean
isomorphism; also, the family $\{s_A(a)\st a\in A\}$ is an open base of $(X,\TT)$. Further, for every Stone space $X$ and   for every $x\in X$, we set
\begin{equation}\label{ux}
u_x=\{P\in CO(X)\st x\in P\}
\end{equation}
(sometimes we will write also $u_x^X$ instead of $u_x$).
Then $u_x\in Ult(CO(X))$ and the map
$$f:X\lra S(CO(X)), \ \ x\mapsto u_x,$$
is a homeomorphism.

When $\unlb=(B,C)$ is a precontact algebra, the pair $(S(B),R_{\unlb})$ is said to be {\em the canonical Stone adjacency space of} $\unlb$.

Now we can obtain the following strengthening of Theorem
\ref{theorem1}:

\begin{theorem}\label{th3}
{\rm (a)} Each precontact algebra $\unlb=(B,C)$ is isomorphic to the canonical
precontact algebra $(CO(X,\TT),C_{R_{\unlb}})$ of the Stone adjacency space
$((X,\TT),R_{\unlb})$, where $(X,\TT)=S(B)$ and for every $u,v\in X$, $uR_{\unlb}v\iff u\times v\sbe C;$
the isomorphism between them is just the Stone map $s_B:B\longrightarrow
CO(X,\TT)$.
Moreover, the relation $C$ satisfies the axiom (Cref) (resp.,
(Csym); (Ctr)) iff the relation $R_{\unlb}$ is reflexive (resp.,
symmetric; transitive).

\smallskip

\noindent{\rm (b)} There exists a bijective correspondence between the class of
all, up to PCA-isomorphism, precontact algebras and the class of
all, up to TAS-iso\-mor\-ph\-ism, Stone adjacency spaces $(X,R)$; namely, for each precontact algebra $\unlb=(B,C)$,
the PCA-isomorphism class $[\unlb]$ of $\unlb$ corresponds to the TAS-isomorphism class of the canonical Stone adjacency space $(S(B), R_{\unlb})$ of $\unlb$,
and for each Stone adjacency space $(X,R)$, the TAS-isomorphism class $[(X,R)]$ of $(X,R)$ corresponds to the PCA-isomorphism class of the canonical precontact algebra $(CO(X),C_R)$ of $(X,R)$ (see (\ref{CR}) for $C_R$).
\end{theorem}

\doc (a) Let $\unlb=(B,C)$  be a precontact algebra, $(X,\TT)=S(B)$ and for every $u,v\in X$, $uR_{\unlb}v\iff u\times v\sbe C$. We will show that $R_{\unlb}$ is a closed relation.  Let $(u,v)\nin R_{\unlb}$. Then there exist $a\in u$ and $b\in v$ such that $a(-C)b$. Hence $u\in s_B(a)$ and $v\in s_B(b)$. Also, $(s_B(a)\times s_B(b))\cap R_{\unlb}=\ems$. Indeed, let $u\ap\in s_B(a)$ and $v\ap\in s_B(b)$; then $a\in u\ap$ and $b\in v\ap$; since $a(-C)b$, we get that $u\ap(-R_{\unlb})v\ap$, i.e. $(u\ap,v\ap)\nin R_{\unlb}$. Therefore, $R_{\unlb}$ is a closed relation. Thus $((X,\TT),R_{\unlb})$ is a Stone adjacency space. We have, by the Stone Representation Theorem, that $CO(X,\TT)=s_B(B)$. Further, we have that for every $a,b\in B$, $s_B(a)C_{R_{\unlb}}s_B(b)\iff (\ex u\in s_B(a))(\ex v\in s_B(b))(u R_{\unlb} v)$. It is easy to see that $(CO(X,\TT),C_{R_{\unlb}})$ is a precontact algebra. We will show that $s_B:(B,C)\longrightarrow
(CO(X,\TT),C_{R_{\unlb}})$ is a PCA-isomorphism. We know, by the Stone Representation Theorem, that $s_B:B\longrightarrow
CO(X,\TT)$ is a Boolean isomorphism. Let $a,b\in B$. Then, using Lemma \ref{lmpc}(a), we obtain that $s_B(a)C_{R_{\unlb}}s_B(b)\iff (\ex u\in s_B(a))(\ex v\in s_B(b))(u R_{\unlb} v)\iff (\ex u,v\in X)((a\in u)\we(b\in v)\we(u R_{\unlb}v))\iff aCb$. Therefore,  $s_B:(B,C)\longrightarrow
(CO(X,\TT),C_{R_{\unlb}})$ is a PCA-isomorphism. The rest follows from Lemma \ref{lmpc}(d,e,f).

\medskip

\noindent(b) Let us set, for every precontact algebra $\unlb=(B,C)$,
$$\Phi(\unlb)=(S(B),R_{\unlb}).$$
Then, by (a),  $\Phi(\unlb)$ is a Stone adjacency space. Further, for every Stone adjacency space $(X,R)$, we set
$$\Psi(X,R)=(CO(X),C_R),$$
where, for every $F,G\in CO(X)$, $FC_R G\iff(\ex x\in F)(\ex y\in G)(xRy)$. Clearly, $\Psi(X,R)$ is a precontact algebra.

Let $\unlb=(B,C)$ be a precontact algebra. Then $\unlb$ is PCA-isomorphic to the precontact algebra $\Psi(\Phi(\unlb))$. Indeed, we have that $\Psi(\Phi(\unlb))=\Psi(S(B),R_{\unlb})=(CO(S(B)),C_{R_{\unlb}})$. Then, by (a), $s_B:(B,C)\longrightarrow \Psi(\Phi(\unlb))$ is a PCA-isomorphism.

Let $(X,R)$ be a Stone adjacency space. Then $(X,R)$ is TAS-isomorphic to the Stone adjacency space $\Phi(\Psi(X,R))$. Indeed, let $B=CO(X)$ and $\unlb=(B,C_R)$. Then
$\Phi(\Psi(X,R))=\Phi(\unlb)=(S(B),R_{\unlb})$. By the Stone Representation Theorem, we have that the map
\begin{equation}\label{stonefux}
f:X\lra S(B), \ x\mapsto u_x,\ \mbox{ is a homeomorphism.}
 \end{equation}
Let $x,y\in X$ and $xR y$. Since, for every $F\in u_x$ and every $G\in u_y$, we have that $x\in F$ and $y\in G$, we get that $u_x R_{\unlb} u_y$, i.e. $f(x) R_{\unlb} f(y)$. Conversely, let $x,y\in X$  and $f(x) R_{\unlb} f(y)$, i.e. $u_x R_{\unlb} u_y$. Suppose that $x(-R)y$. Then $(x,y)\nin R$. Since $R$ is a closed relation, there exist $F,G\in CO(X)$ such that $x\in F$, $y\in G$ and $(F\times G)\cap R=\ems$. Then $F\in u_x$, $G\in u_y$ and if $x\ap\in F$, $y\ap\in G$ then $(x\ap,y\ap)\nin R$, i.e. $x\ap(-R)y\ap$. This implies that $u_x(- R_{\unlb}) u_y$, a contradiction. Therefore, $xRy$. Hence,  $f:(X,R)\lra \Phi(\Psi(X,R))$ is a TAS-isomorphism.
\sqs

\begin{nota}\label{notaclr}
\rm
Let $B$ be a Boolean algebra and $(X,\TT)=S(B)$. Then we denote by $\PCR(B)$ (resp., $\CR(B)$) the set  of all precontact (resp., contact) relations on $B$  and by
$\CLR(X,\TT)$ (resp., $\CLRSR(X,\TT)$) the set of all closed relations  (resp.,  all reflexive and symmetric closed relations) on $(X,\TT)$.
\end{nota}

\begin{cor}\label{contrelbij}
Let $B$ be a Boolean algebra and $(X,\TT)=S(B)$. Then there exists an isomorphism between the ordered  sets $(\PCR(B),\sbe)$ (resp., $(\CR(B),\sbe)$) and  $(\CLR(X,\TT),\sbe)$ (resp., $(\CLRSR(X,\TT),\sbe)$).
\end{cor}

\doc By  Theorem \ref{th3}(a), for every precontact relation $C$ on $B$, the relation $R_{\unlb}$ (defined there), where $\unlb=(B,C)$, is a closed relation on $(X,\TT)$. Now we define
the correspondences $$\Phi\ap:\PCR(B)\lra\CLR(X,\TT), \ \ C\mapsto R_{\unlb},\ \ \mbox{ and }$$  $$\Psi\ap:\CLR(X,\TT)\lra \PCR(B), \ \
R\mapsto C^R,$$ where, for any $a,b\in B$,
$$aC^R b\ \ \iff \ \ s_B(a)\, C_R \, s_B(b).$$ Then, using   Theorem \ref{th3}(a), we get that $\Psi\ap\circ\Phi\ap= id_{\PCR(B)}$. Further, arguing as in the proof of Theorem \ref{th3}(b), we will show that $\Phi\ap\circ\Psi\ap= id_{\CLR(X,\TT)}$. Let $R\in\CLR(X,\TT)$. Then $\Phi\ap(\Psi\ap(R))=R_{(B,C^R)}$. Obviously, for every $u,v\in X$, we have that $u R_{(B,C^R)} v\iff u\times v\sbe C^R\iff (\fa a\in u)(\fa b\in v)(a C^R b)\iff (\fa a\in u)(\fa b\in v)(s_B(a)\, C_R \, s_B(b))$. Since $s_B(a)\, C_R\, s_B(b)\iff (\ex u\ap\in s_B(a))(\ex v\ap\in s_B(b))(u\ap R v\ap)$, we get immediately   that $R\sbe R_{(B,C^R)}$. Let now $u,v\in X$ and $u R_{(B,C^R)} v$. Suppose that $u(-R)v$. Since $R$ is a closed relation, there exist $a,b\in B$ such that $(u,v)\in s_B(a)\times s_B(b)\sbe X^2\stm R$. Then $a\in u$, $b\in v$ and $(s_B(a)\times s_B(b)) \cap R=\ems$. Thus, for every $u\ap\in s_B(a)$ and for every $v\ap\in s_B(b)$, we get that $u\ap(-R) v\ap$. Hence $s_B(a)(- C_R) s_B(b)$,  a contradiction. So, $uRv$. Therefore, $\Phi\ap\circ\Psi\ap= id_{\CLR(X,\TT)}$.
Hence $\Phi\ap$ and $\Psi\ap$ are bijections.
Note that, by Theorem \ref{th3}(a), $C$ is a contact relation on $B$ iff $R_{(B,C)}$ is, in addition, a reflexive and symmetric relation on $X$. Finally, it is clear from the corresponding definitions that for any $C_1,C_2\in\PCR(B)$, $C_1\sbe C_2$ iff $R_{(B,C_1)}\sbe R_{(B,C_2)}$ iff $\Phi\ap(C_1)\sbe\Phi\ap(C_2)$.
\sqs

As it is shown in \cite{DV}, there is no bijective correspondence
between the classes of all, up to corresponding
isomorphisms, precontact algebras and adjacency spaces. Hence, the role of the topology in Theorem
\ref{th3} is essential. However, Theorem \ref{th3} is not
completely satisfactory because the
 representation of the precontact algebras $(B,C)$ obtained here does not give a
topological representation of the  contact algebras $(B,C^\#)$
generated by $(B,C)$; we would like to have an isomorphism $f$
such that, for every $a,b\in B$, $aCb$ iff $f(a)C_R f(b)$, and
$aC^\# b$ iff $f(a)\cap f(b)\not =\emptyset$ (see (\ref{CR}) for
$C_R$). The isomorphism $s_B$ in Theorem \ref{th3} is not of this
type. Indeed, there are many examples of contact algebras $(B,C)$
where $a.b=0$ (and hence $s_B(a)\cap s_B(b)=\emptyset$) but $aCb$
(note that $C$ and $C^{\#}$ coincide for contact algebras). We
 now  construct some natural topological objects which
correspond bijectively to the precontact algebras and satisfy the
above requirement. In the case when $(B,C)$ is a contact algebra,
 we will show that these
topological objects are just topological pairs
satisfying some natural conditions. In such a way we will
obtain new representation theorems for the contact algebras,
 completely different from
those given in \cite{VDDB,DV1,DV2,DW}.

\section{2-Precontact spaces}

Let us start with recalling the following well known statement  (see, e.g., \cite{CNG},
p.271).

\begin{lm}\label{isombool}
Let $X$ be a dense subspace of a topological space $Y$. Then the
functions $$r:RC(Y)\lra RC(X),\  F\mapsto F\cap X,$$ and
$$e:RC(X)\lra
RC(Y),\  G\mapsto \cl_Y(G),$$ are Boolean isomorphisms between Boolean
algebras $RC(X)$ and $RC(Y)$, and $e\circ r=id_{RC(Y)}$, $r\circ
e=id_{RC(X)}$. (We will sometimes write $r_{X,Y}$ (resp., $e_{X,Y}$) instead of $r$ (resp., $e$).)
\end{lm}

\begin{defi}\label{toppair}
\rm
(a) Let $X$ be a topological space and $X_{0}$ be  a dense
subspace of $X$. Then the pair
 $(X,X_{0})$ is called a\/ {\em topological pair}.

\medskip

\noindent(b) Let $(X,X_0)$ be a topological pair. Then we set
\begin{equation}\label{rcx0}
  RC(X,X_{0})=\{cl_X(A)\ |\ A\in CO(X_0)\}.
  \end{equation}
\end{defi}

\begin{lm}\label{lmrcx0}  Let $(X,X_0)$ be a topological pair.
Then $RC(X,X_0)\subseteq RC(X)$; the set $RC(X,X_{0})$ with the
standard Boolean operations on the regular closed subsets of $X$
is a  Boolean subalgebra of  $RC(X)$; $RC(X,X_0)$ is isomorphic to
the Boolean algebra $CO(X_{0})$; the sets $RC(X)$ and
$RC(X,X_0)$  coincide iff $X_0$ is an extremally disconnected
space. If $$C_{(X,X_0)}$$ is the restriction of the contact relation
$C_{X}$ (see Lemma \ref{rcl}) to $RC(X,X_{0})$, then $(RC(X,X_{0}),
C_{(X,X_{0})})$ is a  contact subalgebra of $(RC(X), C_{X})$.
\end{lm}

\doc
Since $CO(X_0)$ is a Boolean subalgebra of the Boolean algebra $RC(X_0)$ and (in the notation of Lemma \ref{isombool}) $e(CO(X_0))=RC(X,X_0)$ and $e$ is a Boolean isomorphism, we get that $RC(X,X_0)\subseteq RC(X)$ and the set $RC(X,X_{0})$ with the
standard Boolean operations on the regular closed subsets of $X$
is a  Boolean subalgebra of  $RC(X)$. Clearly, the restriction $e_0=e|_{CO(X_0)}:CO(X_0)\lra RC(X,X_0)$ is a Boolean isomorphism. Using the above arguments, we get
that  $RC(X)=RC(X,X_0)$  iff $RC(X_0)=CO(X_0)$. As it is well known,  the later equality is true iff
$X_0$ is an extremally disconnected
space; hence $RC(X)=RC(X,X_0)$  iff $X_0$ is an extremally disconnected
space.
Finally, it is obvious that $(RC(X,X_{0}),
C_{(X,X_{0})})$ is a  contact subalgebra of $(RC(X), C_{X})$.
\sqs

\begin{nota}\label{sigx}
\rm
Let $(X,\tau)$ be a topological space, $X_0$ be a subspace of
$X$, $x\in X$ and $B$ be a subalgebra of the Boolean algebra $(RC(X),+,.,*,\ems,X)$ defined in \ref{rcdef}. We put
\begin{equation}\label{sg}
\sigma_x^B=\{F\in B\ |\ x\in F\};\ \
\Gamma_{x,X_0}=\{F\in
CO(X_0)\ |\ x\in cl_X(F)\}.
\end{equation}

We set also

\begin{equation}\label{sg1}
\nu_x^B=\{F\in B\st x\in\int_X(F)\}.
\end{equation}

When $B=RC(X)$, we will often write simply $\s_x$ and $\nu_x$ instead of, respectively, $\s_x^B$ and $\nu_x^B$; in this case we will sometimes use the notation $\s_x^X$ and $\nu_x^X$ as well.
\end{nota}

\begin{defi}\label{precontactspace} \textsc{(2-Precontact spaces.)}
\rm

\noindent(a) A triple $\unlx=(X,X_0,R)$ is  called a {\em
2-precontact space} (abbreviated as PCS) if the following
conditions are satisfied:
\begin{description}
\item[$(PCS 1)$]\quad  $(X,X_{0})$ is a topological pair and
 $X$ is a  $T_{0}$-space;

\item[$(PCS 2)$]\quad  $(X_{0},R)$ is a Stone adjacency space;

\item[$(PCS 3)$]\quad $RC(X,X_0)$ is a closed base for $X$;

\item[$(PCS 4)$]\quad For every $F,G\in CO(X_0)$, $\cl_X(F)\cap
\cl_X(G)\neq\emptyset$ implies that $F(C_R)^\# G$ (see (\ref{CR})
for $C_R$);

\item[$(PCS 5)$]\quad If $\Gamma\in Clans(CO(X_0),C_R)$ then
there exists a point $x\in X$ such that $\Gamma = \Gamma_{x,X_0}$
(see (\ref{sg}) for $\Gamma_{x,X_0}$).
\end{description}

\noindent(b) Let\/ $\unlx=(X,X_0,R)$ be a 2-precontact space. Define,
for every $F,G\in RC(X,X_0)$,
$$F\ C_{\unlx}\ G \iff ((\ex x\in
F\cap X_0) (\ex y\in G\cap X_0) (xRy)).$$
 Then the
precontact algebra
$$\unlb(\unlx)= (RC(X,X_{0}), C_{\unlx})$$
 is
said to be the\/ {\em canonical precontact algebra of\/
$\unlx$}.

\smallskip

\noindent(c) A 2-precontact space\/ $\unlx=(X,X_0,R)$ is called\/ {\em
reflexive} (resp.,\/ {\em symmetric}; {\em transitive}) if the
relation $R$ is reflexive (resp., symmetric; transitive);\/
$\unlx$ is called\/ {\em connected} if the space $X$ is
connected.

\smallskip

\noindent(d)  Let\/ $\unlx=(X,X_{0},R)$ and\/
$\widehat{\unlx}=(\widehat{X},\widehat{X}_{0},\widehat{R})$
be two 2-precontact spaces. We say that\/ $\unlx$ and\/
$\widehat{\unlx}$ are\/ {\em PCS-isomorphic} (or, simply,\/
{\em isomorphic}) if there exists a homeomorphism
$f:X\longrightarrow\widehat{X}$ such that:
\begin{description}
\item[{\rm (ISO1)}]\quad   $f(X_{0})=\widehat{X}_{0}$; and

\item[{\rm (ISO2)}]\quad $(\forall x,y\in X_{0})(xRy\leftrightarrow
f(x)\widehat{R} f(y))$.
\end{description}
\end{defi}

\begin{rem}\label{rem2pre}
\rm
It is very easy to see that the canonical precontact algebra of a 2-precontact space, defined in Definition \ref{precontactspace}(b),  is indeed a precontact algebra.
\end{rem}

\begin{pro}\label{proposition1}
{\rm (a)} Let $(X,X_0,R)$ be a
2-precontact space. Then $X$ is a semire\-gular space and, for
every $F,G\in CO(X_0)$,
\begin{equation}\label{clfg}
cl_X(F)\cap cl_X(G)\neq\emptyset \mbox{ iff }F(C_{R})^\# G.
\end{equation}

\smallskip

\noindent(b) Let\/ $\unlx=(X,X_0,R)$ and\/ $\widehat{\unlx}=(\widehat{X},\widehat{X}_{0},\widehat{R})$ be two
isomorphic 2-precontact spaces. Then the corresponding canonical
precontact algebras\/ $\unlb(\unlx)$  and\/
$\unlb(\widehat{\unlx})$ are PCA-isomorphic.
\end{pro}

\doc (a) By the axiom (PCS3), the family $RC(X, X_0)$ is a closed base for the space $X$. Since $RC(X,X_0)\sbe RC(X)$ (see Lemma \ref{lmrcx0}(a)), we get that $X$ is a semiregular space.

Let $F,G\in CO(X_0)$ and $F(C_R)^\# G$.  The pair $(CO(X_0),C_R)$ is a precontact algebra (see Proposition \ref{prop1}(a)). Hence, by Proposition \ref{lmpc}(c), there exists a clan $\GA$ in $(CO(X_0),C_R)$ such that $F,G\in\GA$. The axiom (PCS5) implies that there exists $x\in X$ such that $\GA=\GA_{x,X_0}$. Thus $x\in cl_X(F)\cap cl_X(G)$. Therefore, $cl_X(F)\cap cl_X(G)\neq\emptyset$. The converse implication follows from the axiom (PCS4).

\medskip

\noindent(b)  This is obvious.
\sqs

\begin{lm}\label{correspondencelemma} \textsc{(Correspondence Lemma.)}
 Let\/ $\unlx=(X,X_{0},R)$ be a 2-precontact space and let
\/ $\unlb(\unlx)= (RC(X,X_{0}), C_{\unlx})$ be the canonical
precontact algebra of\/ $\unlx$. Then the following
equivalences hold:
\begin{description}

   \item[{\rm (a)}]\quad The space\/ $\unlx$ is reflexive iff
    the algebra\/  $\unlb(\unlx)$
satisfies the axiom  $(Cref)$;

\item[{\rm (b)}]\quad The space\/  $\unlx$ is symmetric iff\/
$\unlb(\unlx)$ satisfies the axiom $(Csym)$.;

\item[{\rm (c)}]\quad The space\/  $\unlx$ is transitive iff\/
$\unlb(\unlx)$ satisfies the axiom $(Ctr)$;

\item[{\rm (d)}]\quad The space\/  $\unlx$ is connected iff\/  $\unlb(\unlx)$
is connected.
\end{description}
\end{lm}

\doc By (the proof of) Lemma \ref{lmrcx0}(a), the map $$\p:CO(X_0)\lra RC(X,X_0),\ \ A\mapsto\cl_X(A),$$ is a Boolean isomorphism. From the definitions of the relations
$C_{\unlx}$ and $C_R$ (see, respectively, Definition \ref{precontactspace}(b) and (\ref{CR})) it follows immediately that the map
\begin{equation}\label{corcx0}
\p:(CO(X_0),C_R)\lra (RC(X,X_0),C_{\unlx})\ \mbox{ is a PCA-isomorphism.}
 \end{equation}
 Now the assertions (a), (b) and (c) follow from Definition \ref{precontactspace}(c), Theorem \ref{th3} and Lemma \ref{lmpc}(d,e,f). Indeed, by the axiom (PCS2), $(X_0,R)$ is a Stone adjacency space; hence, by Theorem \ref{th3}(b) and in the notation of its proof, $\Psi(X_0,R)=(CO(X_0),C_R)$ and $\Phi(\Psi(X_0,R))=(S(CO(X_0)),R_{(CO(X_0),C_R)})$ is TAS-isomorphic to the  Stone  adjacency space $(X_0,R)$; now we can apply Lemma \ref{lmpc}(d,e,f) to the precontact algebra $(CO(X_0),C_R)$ and its canonical adjacency space $(Ult(CO(X_0)),R_{(CO(X_0),C_R)})=(S(CO(X_0)),R_{(CO(X_0),C_R)})$.
So, the assertions (a), (b) and (c) are proved.
Let us prove the assertion (d). Let $X$ be a connected space. Suppose that the precontact algebra $(CO(X_0),C_R)$ is not connected. Then there exists an $F\in CO(X_0)$ such that $F\neq\ems,X_0$, $F(-C_R)F^*$ and $F^*(-C_R)F$. Since $F.F^*=\ems=0$, we get that $F(-(C_R)^\# )F^*$. Thus, by (\ref{clfg}), $\cl_X(F)\cap\cl_X(F^*)=\ems$. Since $F^*=X_0\stm F$ and $X_0$ is dense in $X$, we obtain that $\cl_X(F)\cup\cl_X(F^*)=\cl_X(F)\cup\cl_X(X_0\stm F)=\cl_X(X_0)=X$. Hence $\cl_X(F)$ is a clopen  subset of $X$. Since $\cl_X(F)\neq\ems,X$, we get a contradiction. Thus, the precontact algebra $(CO(X_0),C_R)$ is  connected. Conversely, let the precontact algebra $(CO(X_0),C_R)$ be  connected. Suppose that $X$ is not connected. Then there exists a clopen in $X$ subset $G$ of $X$ such that $G\neq\ems,X$. Let $F=X_0\cap G$. Then $F\in CO(X_0)$ and $F\neq\ems,X_0$. Thus, $FC_RF^*$ or $F^*C_RF$. Hence $F(C_R)^\# F^*$. Using (\ref{clfg}), we get that $\cl_X(F)\cap\cl_X(F^*)\nes$. Since $\cl_X(F)=\cl_X(G\cap X_0)=\cl_X(G)=G$ and, analogously, $\cl_X(F^*)=\cl_X(X_0\stm F)=\cl_X((X\stm G)\cap X_0)=X\stm G$, we get a contradiction. Hence, $X$ is connected.
\sqs

\section{The canonical 2-precontact space of a precontact algebra}

In this section we will associate with each precontact algebra a
2-precontact space.

\begin{defi}\label{canspa}
\rm
 Let $\unlb=(B,C)$ be a precontact algebra. We associate
 with
 $\unlb$ a 2-precontact space
$$\unlx(\unlb)=(X,X_{0},R),$$
called the\/ {\em canonical
2-precontact space of\/ $\unlb$}, as follows:
\begin{description}

\item[$\bullet$] \quad  $X=Clans(\unlb)$  and $X_{0}=Ult(B)$;

\item[$\bullet$] \quad    The topology $\tau$ on the set  $X$ is defined in
the following way:  the family $$\{g_B(a)\
| \ a\in B\},$$ where,
for any $a\in B$,
\begin{equation}\label{fg}
g_B(a)=\{ \Gamma\in X \ | \ a\in \Gamma\},
\end{equation}
is a closed base of $\tau$. The topology on $X_0$ is the subspace topology induced by $(X,\tau)$.

\item[$\bullet$] \quad  $R=R_{\unlb}$  (see (\ref{ultra}) for the notation $R_{\unlb}$).
\end{description}
\end{defi}

\begin{rem}\label{rem2prc41}
\rm
Note that, in the notation of Definition \ref{canspa}, setting,
for every $a\in B$, $$g_{0}^B(a)=g_B(a)\cap X_0,$$
we get that the family
$\{g_0^B(a)\ | \ a\in B\}$ is a closed base of $X_0$ and
$g_0^B(a)=s_B(a)$, where $s_B:B\longrightarrow CO(X_0)$ is the Stone
map.
\end{rem}

\begin{pro}\label{proposition2}  Let\/ $\unlb=(B,C)$ be a precontact algebra.
 Then the  canonical
2-precontact space \/ $\unlx(\unlb)=(X,X_{0},R)$ of\/ $\unlb$
  defined above is indeed a 2-precontact
space.
\end{pro}

\doc Since $Clans(B,C)\equiv Clans(B,C^\#)$ and $(B,C^\#)$ is a contact algebra, we can use \cite[Lemma 5.1(i), Lemma 5.7(i) and Lemma 5.3(ii)]{DV1} which imply that the family $\{g_B(a)\ | \ a\in B\}$ can be taken as a closed base for a topology on the set $X$, that $(X,\tau)$ is a semiregular compact $T_0$-space and that, for every $a\in B$, $g_B(a)\in RC(X)$.

Let us set, for every $a\in B$, $h_B(a)=X\stm g_B(a)$. Then $$h_B(a)=\{\GA\in X\st a\nin \GA\}$$
and $\{h_B(a)\st a\in B\}$ is an open base of $(X,\tau)$.
   Let us show that $X_0$ is dense in $(X,\tau)$. Let $a\in B$ and $a\neq 1$ (i.e., $h_B(a)\nes$). Then $a^*\neq 0$ and thus there exists an ultrafilter $u$ in $B$ such that $a^*\in u$. Thus $a\nin u$ and hence $u\in h_B(a)$. Therefore, $h_B(a)\cap X_0\nes$. Hence, $X_0$ is dense in $(X,\tau)$. So, the axiom (PCS1) is satisfied.
   Since $g_0^B(a)=s_B(a)$, for every $a\in B$ (see Remark \ref{rem2prc41}),  Theorem \ref{th3} implies  that $((X_0,\tau|_{X_0}),R)$ is a Stone adjacency space. Thus, the axiom (PCS2) is also satisfied. Further, we have that for every $a\in B$, $g_B(a)\in RC(X)$; thus, using Lemma \ref{isombool}, we get that
\begin{equation}\label{g0acl}
g_B(a)=\cl_X(g_0^B(a)).
 \end{equation}
 Since $CO(X_0)=\{g_0^B(a)\st a\in B\}$, we obtain  that $RC(X,X_0)=\{g_B(a)\st a\in B\}$. Hence, the axiom (PCS3) is satisfied.

Let $F,G\in CO(X_0)$. Then there exist $a,b\in B$ such that $F=s_B(a)(=g_0^B(a))$ and $G=s_B(b)(=g_0^B(b))$. Let $\cl_X(F)\cap\cl_X(G)\nes$. Then, by (\ref{g0acl}), $g_B(a)\cap g_B(b)\nes$. Hence, there exists $\GA\in g_B(a)\cap g_B(b)$. Then $a,b\in\GA$ and thus, $aC^\# b$. There are three possibilities:

\smallskip

\noindent 1) $aCb$: then, by Lemma \ref{lmpc}(a), there exist $u,v\in X_0$ such that $a\in u$, $b\in v$ and $uRv$; since $u\in F$ and $v\in G$, we get that $FC_RG$ and therefore, $F(C_R)^\#G$;

\smallskip

\noindent 2) $bCa$: then, by Lemma \ref{lmpc}(a), there exist $u,v\in X_0$ such that $b\in u$, $a\in v$ and $uRv$; since $v\in F$ and $u\in G$, we get that $GC_RF$ and therefore, $F(C_R)^\#G$;

 \noindent 3)  $a.b\neq 0$: then, as it is well known (see, e.g., \cite[Corollary 2.17]{kop89}), there exists $u\in X_0$ such that $a.b\in u$; thus $a,b\in u$, i.e., $u\in F\cap G$ and therefore, $F(C_R)^\#G$.

So, the axiom (PCS4) is satisfied.

Let $\GA\in Clans(CO(X_0),C_R)$. We have, by Theorem \ref{th3}(a), that
\begin{equation}\label{scpcadiez}
s_{(B,C)}:(B,C)\lra (CO(X_0),C_R),\ \ b\mapsto s_B(b),\ \mbox{ is a PCA-isomorphism.}
\end{equation}
 Set $x=s_{(B,C)}\inv(\GA)$. Then $x$ is a clan of  $(B,C)$, i.e. $x\in X$. We will show that $\GA=\GA_{x,X_0}$ (see (\ref{sg}) for the notation $\GA_{x,X_0}$). We have to prove that $(\fa F\in CO(X_0))((x\in \cl_X(F))\iff (F\in\GA))$. Let $F\in CO(X_0)$. Then there exists $a\in B$ such that $F=s_{(B,C)}(a)$. Let $F\in\GA$. Then $s_{(B,C)}(a)\in\GA$ and thus $s_{(B,C)}\inv(s_{(B,C)}(a))\in s_{(B,C)}\inv(\GA)$, i.e. $a\in x$ and hence $x\in g_B(a)=\cl_X(s_B(a))=\cl_X(F)$. Conversely, let $x\in\cl_X(F)$. Since $g_B(a)=\cl_X(s_B(a))$, we get that $a\in x$. Then $s_B(a)\in s_{(B,C)}(x)=\GA$. Therefore, $F\in\GA$. Hence there exists $x\in X$ such that $\GA=\GA_{x,X_0}$. Thus the axiom (PCS5) is satisfied.
\sqs

\begin{pro}\label{proposition3}  Let\/ $\unl{B_1}$
and\/ $\unl{B_{2}}$ be two
 isomorphic precontact algebras. Then the corresponding canonical 2-precontact
 spaces\/ $\unlx(\unl{B_1})$ and\/ $\unlx(\unl{B_2})$ are isomorphic.
\end{pro}

\doc It is obvious.
\sqs

\begin{pro}\label{proposition4} Let $ \unlb$ be a precontact
algebra, let $\unlx(\unlb)$ be the canonical 2-precon\-t\-act
space of  $\unlb$ and let    $\unl{B'}$ be the
canonical precontact algebra of the 2-precontact space
$\unlx(\unlb)$. Then the contact algebras  $\unlb$ and\/  $\unl{B'}$ are
PCA-isomorphic.
\end{pro}

\doc It follows from (\ref{scpcadiez}) and (\ref{corcx0}).
\sqs

\begin{lm}\label{conect}   \textsc{(Topological characterization of
the connectedness of} $\unlb$.) Let $\unlb=(B,C)$ be a
precontact algebra and $\unlx(\unlb)$ be the canonical
2-precontact space of $\unlb$. Then $\unlb$ is
connected iff $\unlx(\unlb)$ is  connected.
\end{lm}

\doc Since $Clans(B,C)\equiv Clans(B,C^\#)$ and $(B,C^\#)$ is a contact algebra, we can use \cite[Lemma 5.7(i3)]{DV1} which implies that
the 2-precontact space $\unlx(\unlb)$ is connected iff the contact algebra  $(B,C^\#)$ is connected. Obviously, $(B,C^\#)$ is connected iff $\unlb$ is connected.
Thus, our assertion is proved.
\sqs

\section{The Main Theorem}

\begin{theorem}\label{mth} \textsc{(Representation theorem for precontact algebras.)}
\begin{description}
\item[{\rm (a)}]\quad  Let $\unlb=(B,C)$ be a precontact algebra
and let
$\unlx(\unlb)=(X,X_0,R)$
 be the canonical 2-precontact
space of $\unlb$. Then the function $g_B:(B,C)\lra 2^X$, defined in
(\ref{fg}), is a PCA-isomorphism  from $(B,C)$ onto the canonical
precontact algebra $(RC(X,X_0),C_{\unlx(\unlb)})$ of $\unlx(\unlb)$. The
same function $g_B$ is a PCA-isomorphism between  contact
algebras $(B,C^\#)$ and $(RC(X,X_0),C_{(X,X_0)})$ (see Lemma \ref{lmrcx0}(a) for $C_{(X,X_0)}$). The sets
 $RC(X)$ and  $RC(X,X_0)$ coincide iff the precontact algebra
$\unlb$ is complete. The algebra $\unlb$ satisfies the
axiom $(Cref)$ (resp., $(Csym); (Ctr))$ iff the 2-precontact
space $\unlx(\unlb)$ is reflexive (resp., symmetric;
transitive). The algebra $\unlb$ is connected iff
$\unlx(\unlb)$ is connected.

\item[{\rm (b)}]\quad There exists a bijective correspondence  between the class of
all, up to PCA-isomor\-ph\-ism, (connected) precontact algebras
and  the class of all, up to PCS-isomorphism, (connected) 2-precontact spaces; namely, for every precontact algebra $\unlb$, the PCA-isomorphism class $[\unlb]$ of $\unlb$ corresponds to the PCS-isomor\-phism class $[\unlx(\unlb)]$ of the canonical 2-precontact space $\unlx(\unlb)$ of $\unlb$, and for every 2-precontact space $\unlx$, the PCS-isomorphism class $[\unlx]$ of $\unlx$ corresponds to the PCA-isomorphism class $[\unlb(\unlx)]$ of the canonical precontact algebra $\unlb(\unlx)$ of $\unlx$.
\end{description}
\end{theorem}

\doc (a) Using (\ref{g0acl}), we get that $\p\circ s_B=g_B$ (see  (\ref{corcx0}) for $\p$ and (\ref{stonemap}) for $s_B$).  Now, we apply Proposition \ref{proposition4} for obtaining that
\begin{equation}\label{giso}
g_B:(B,C)\lra (RC(X,X_0),C_{\unlx(\unlb)})\ \mbox{ is a PCA-isomorphism.}
\end{equation}
Further, using (\ref{clfg}), we get that
\begin{equation}\label{gisoc}
g_B:(B,C^\#)\lra (RC(X,X_0),C_{(X,X_0)})\ \mbox{ is a CA-isomorphism.}
\end{equation}

By Lemma \ref{lmrcx0}(a), we have that the sets
 $RC(X)$ and  $RC(X,X_0)$ coincide iff $X_0$ is an extremally disconnected
space. Since $X_0=S(B)$, we have (by, e.g., \cite[Proposition 7.21]{kop89}) that $X_0$ is an extremally disconnected
space iff $B$ is a complete Boolean algebra.
Hence, the sets
 $RC(X)$ and  $RC(X,X_0)$ coincide iff the precontact algebra
$\unlb$ is complete. All other assertion in (a) follow from Lemma \ref{correspondencelemma} and (\ref{giso}).

\medskip

\noindent(b) Let us denote by $\PCA$ the set of all, up to PCA-isomor\-ph\-ism,  precontact algebras and by $\PCS$ the set of
all, up to PCS-isomorphism,  2-precontact spaces. We will define two correspondences
$$\Phi_2:\PCA\lra \PCS\ \ \mbox{ and }\ \ \Psi_2:\PCS\lra \PCA$$
and we will show that their compositions $\Phi_2\circ\Psi_2$ and $\Psi_2\circ\Phi_2$ are equal to the corresponding identities.
We set, for every precontact algebra $\unlb=(B,C)$, $$\Phi_2([\unlb])=[\unlx(\unlb)],$$ where $\unlx(\unlb)$ is the canonical 2-precontact space of $\unlb$ (see
Definition \ref{canspa}), $[\unlb]$ is the class of all precontact algebras which are PCA-isomorphic to the precontact algebra $\unlb$, and, analogously,
$[\unlx(\unlb)]$ is the class of all 2-precontact spaces which are PCS-isomorphic to the 2-precontact space $\unlx(\unlb)$. Further, for every 2-precontact space $\unlx=(X,X_0,R)$, we set $$\Psi_2([\unlx])=[\unlb(\unlx)],$$ where $\unlb(\unlx)$ is the canonical precontact algebra of $\unlx$ (see Definition \ref{precontactspace}(b)).
It is easy to see that the correspondences $\Phi_2$ and $\Psi_2$ are well-defined.

Using (\ref{giso}), we get that for every precontact algebra $\unlb=(B,C)$, $\Psi_2(\Phi_2([\unlb]))=[\unlb]$. Thus we obtain that $\Psi_2\circ\Phi_2=id_{\PCA}$.

 We will now prove that $\Phi_2\circ\Psi_2=id_{\PCS}$. Let $\unlx=(X,X_0,R)$ be a 2-precontact space. Set $(B,C)=(CO(X_0),C_R)$; then $(B,C)$ is PCA-isomorphic to the canonical precontact algebra of $\unlx$ (see Definition \ref{precontactspace}(b)). Let $(\wh{X},\wh{X}_0,\wh{R})$ be the canonical 2-precontact space of $(B,C)$ (see
Definition \ref{canspa} and (\ref{corcx0})). Then $\wh{X}=Clans(B,C)$, $\wh{X}_0=Ult(B)$ and $\wh{R}=R_{(B,C)}$. For every $x\in X$, set
$$f(x)=\{a\in B\st x\in\cl_X(a)\}.$$
Then $f(x)\nes$ (by (PCS3)) and $f(x)$ is a clan in $(B,C)$. Indeed, we have that: 1) $0\nin f(x)$; 2) if $a\in f(x)$, $b\in B$ and $a\le b$, then $x\in\cl_X(a)\sbe\cl_X(b)$, and thus $b\in f(x)$; 3) if $a+b\in f(x)$ then $x\in\cl_X(a\cup b)=\cl_X(a)\cup\cl_X(b)$, and hence $x\in\cl_X(a)$ or $x\in\cl_X(b)$, i.e. $a\in f(x)$ or $b\in f(x)$; 4) if $a,b\in f(x)$ then $x\in\cl_X(a)\cap\cl_X(b)$, and thus, by (PCS4), $aC^\# b$. So, $f(x)\in\wh{X}$. Hence,
$f:X\lra\wh{X}.$ We will show that
$$f:(X,X_0,R)\lra (\wh{X},\wh{X}_0,\wh{R}),\ \ x\mapsto f(x),\ \ \mbox{ is a PCS-isomorphism.}$$

Let $x,y\in X$ and $x\neq y$. Since, by (PCS1), $X$ is a $T_0$-space, there exists an open subset $U$ of $X$ such that $|U\cap\{x,y\}|=1$. Suppose that $x\in U$. Then $y\nin U$. According to (PCS3), there exists $a\in B$ such that $x\in X\stm\cl_X(a)\sbe U$. Thus $y\in\cl_x(a)$ and $x\nin\cl_X(a)$. Therefore, $a\in f(y)\stm f(x)$. Hence $f(x)\neq f(y)$. If $y\in U$, then we argue analogously. So, $f$ is a injection.

Let $\GA\in\wh{X}$. Then, by (PCS5), there exists $x\in X$ such that $\GA=\{a\in B\st x\in\cl_X(a)\}$. Hence $\GA=f(x)$. Thus $f$ is a surjection. So, $f$ is a bijection.

Let $a\in B$. Then $f(\cl_X(a))=\{f(x)\st x\in\cl_X(a)\}=\{f(x)\st a\in f(x)\}$. Since $f$ is a surjection, we get that $f(\cl_X(a))=g_B(a)$ (see (\ref{fg}) for the notation $g_B(a)$). Since $f$ is a bijection, we have also that $f\inv(g_B(a))=\cl_X(a)$. Now, using (PCS3) and the fact that $\{g_B(a)\st a\in B\}$ is a closed base of $\wh{X}$ (see Definition \ref{canspa}), we get that $f$ is a homeomorphism.

For every $x\in X_0$, we have that $f(x)=\{F\in CO(X_0)\st x\in\cl_X(F)\}=\{F\in CO(X_0)\st x\in F\}=u_x\in Ult(B)=\wh{X}_0$ (see (\ref{ux}) for the notation $u_x$). Hence $f(X_0)\sbe\wh{X}_0$. For proving the inverse inclusion, let $u\in \wh{X}_0$. Then $u\in Ult(CO(X_0))$. Now, by (PCS2), there exists $x\in X_0$ such that $x\in \bigcap u$. Then $u\sbe u_x$ and, hence, $u=u_x$. Since, as we have already seen, $u_x=f(x)$, we get that  $f(X_0)\spe\wh{X}_0$. Therefore, $f(X_0)=\wh{X}_0$.

Let $x,y\in X_0$. Then $f(x)=u_x$, $f(y)=u_y$. We have that $u_x,u_y\in\wh{X}_0$ and $u_x\wh{R}u_y\iff u_x\times u_y\sbe C_R$. Hence, ($u_x\wh{R}u_y$) $\iff$ (for every $F,G\in CO(X_0)$ such that $x\in F$ and $y\in G$, there exist $x\ap\in F$ and $y\ap\in G$ with $x\ap R y\ap$). Therefore, if $xRy$ then, obviously, $f(x)\wh{R} f(y)$.
Let now $f(x)\wh{R} f(y)$. Suppose that $x(-R)y$. Then $(x,y)\nin R$. Applying (PCS2), we get that there exist $F,G\in CO(X_0)$ such that $x\in F$, $y\in G$ and $(F\times G)\cap R=\ems$. Thus $F\in u_x$, $G\in u_y$ and for every $x\ap\in F$ and every $y\ap\in G$ we have that $x\ap(-R)y\ap$. This means that $F(-C_R)G$ and, hence,
$u_x(-\wh{R})u_y$, i.e. $f(x)(-\wh{R}) f(y)$, a contradiction. Therefore, $xRy$.

All this shows that $f$ is a PCS-isomorphism. Hence $\Phi_2(\Psi_2([\unlx]))=[\unlx].$ Thus, $\Phi_2\circ\Psi_2=id_{\PCS}$.
Therefore,
\begin{equation}\label{pcapcsbj}
\Phi_2:\PCA\lra \PCS\ \ \mbox{ is a bijection.}
\end{equation}

The statement for connected precontact algebras follows from (\ref{pcapcsbj})  and Lemma \ref{conect}.
\sqs

\begin{cor}\label{cormth}
If $\unlx=(X,X_0,R)$ is a 2-precontact space then $X$ is a compact space.
\end{cor}

\doc By Theorem \ref{mth}, there exists a precontact algebra $\unlb=(B,C)$ such that the 2-precontact space $\unlx$ is isomorphic to the canonical 2-precontact space $\unlx(\unlb)$ of $\unlb$. Then, by \cite[Lemma 5.7(i2)]{DV1}, $X$ is a compact space (see also the proof of Proposition \ref{proposition2}).
\sqs

\section{2-Contact spaces and a new representation theorem for contact algebras}

\begin{pro}\label{subrel}
Let $X_0$ be a subspace of a topological space  $X$. For every
$F,G\in CO(X_0)$, set
\begin{equation}\label{deltarel}
F\delta_{(X,X_0)}G \mbox{ iff }cl_X(F)\cap cl_X(G)\neq\emptyset.
\end{equation}
Then $(CO(X_0),\delta_{(X,X_0)})$ is a contact algebra.
\end{pro}

\doc Clearly, the relation $\delta_{(X,X_0)}$ satisfies the axioms (C0), (Cref) and (Csym) (see Definition \ref{precontact} for them). It is easy to see that it satisfies the axiom (C+) as well. Hence, $(CO(X_0),\delta_{(X,X_0)})$ is a contact algebra.
\sqs

\begin{defi}\label{contactspace}\textsc{(2-Contact spaces.)}
\rm
(a) A topological pair  $(X,X_0)$ is  called a\/ {\em 2-contact
space} (abbreviated as CS) if the following conditions are
satisfied:
\begin{description}
\item[$(CS 1)$]\quad
 $X$ is a  $T_{0}$-space;

\item[$(CS 2)$]\quad  $X_{0}$ is a Stone space;

\item[$(CS 3)$]\quad $RC(X,X_0)$ is a closed base for $X$;

\item[$(CS 4)$]\quad If $\Gamma\in
Clans(CO(X_0),\delta_{(X,X_0)})$ (see (\ref{deltarel}) for the notation $\delta_{(X,X_0)}$) then there exists a point
$x\in X$ such that $\Gamma =\Gamma_{x,X_0}$ (see (\ref{sg}) for
$\Gamma_{x,X_0}$).
\end{description}

\noindent A 2-contact space $(X,X_0)$ is called\/ {\em connected}
if the space $X$ is connected.

\medskip

\noindent(b) Let  $(X,X_0)$ be a 2-contact space. Then the contact algebra
 $${\unlbc(X,X_0)}= (RC(X,X_{0}), C_{(X,X_0)})$$
(see Lemma \ref{lmrcx0}(a) for the notation $C_{(X,X_0)}$)  is said to be
the\/ {\em canonical contact algebra of the 2-contact space
$(X,X_0)$}.

\medskip

\noindent(c) Let $\unlb=(B,C)$ be a contact algebra, $X=Clans(B,C)$, $X_0=Ult(B)$
and $\tau$ be the topology on $X$ described in Definition
\ref{canspa}. Take  the subspace topology on $X_0$. Then the pair
$$\unlxc(\unlb)=(X,X_0)$$
is called the\/ {\em canonical 2-contact space of the
contact algebra} $(B,C)$.

\medskip

\noindent(d)  Let   $(X,X_{0})$ and $(\widehat{X},\widehat{X}_{0})$ be two
2-contact spaces. We say that  $(X,X_0)$ and
$(\widehat{X},\widehat{X}_0)$ are\/ {\em CS-isomorphic} (or,
simply,\/ {\em isomorphic}) if there exists a homeomorphism
$f:X\longrightarrow\widehat{X}$ such that
  $f(X_{0})=\widehat{X}_0$.
\end{defi}

\begin{rem}\label{remcanconalg}
\rm
Note that, by Lemma \ref{lmrcx0}(a), the canonical contact algebra of a 2-contact space is indeed a contact algebra.
\end{rem}

\begin{exa}\label{stspexa}
\rm
Let $X$ be a Stone space. Then the pair $(X,X)$ is a 2-contact space. Indeed, the conditions (CS1)-(CS3) are obviously fulfilled. Let $B=CO(X)$. Obviously, $\d_{(X,X)}=\rho_s^B$ (see Example \ref{extrcr} for $\rho_s^B$). Hence, using \cite[Corollary 3.3]{DV1} or arguing  directly, we get that $Clans(B,\d_{(X,X)})=Ult(B)=\{u_x\st x\in X\}=\{\GA_{x,X}\st x\in X\}$. Thus, the condition (CS4) is also fulfilled.
\end{exa}

\begin{lm}\label{conprecon}
For every 2-contact space $(X,X_0)$ there exists a unique
reflexive and symmetric binary relation $R$ on $X_0$ such that
$(X,X_0,R)$ is a 2-precontact space.
\end{lm}

\doc Let $(X,X_0)$ be a 2-contact space. For every $x,y\in X_0$, set
\begin{equation}\label{xrycont}
xRy\iff ((\fa F\in u_x)(\fa G\in u_y)(\cl_X(F)\cap\cl_X(G)\nes)),
\end{equation}
where  $u_x=\{A\in CO(X_0)\st x\in A\}$ and analogously for $u_y$ (see (\ref{ux})). Clearly, $R$ is a reflexive and symmetric relation on $X_0$. We will show that for every $F,G\in CO(X_0)$,
\begin{equation}\label{clxfclxgnescr}
 \cl_X(F)\cap\cl_X(G)\nes\iff FC_RG
 \end{equation}
(recall that $FC_RG\iff((\ex x\in F)(\ex y\in G)(xRy))$).
Indeed, let $F,G\in CO(X_0)$ and
 let $\cl_X(F)\cap\cl_X(G)\nes$. Then $F\d_{(X,X_0)}G$. Now, using Proposition \ref{subrel} and Lemma \ref{lmpc}(c), we get that there exists a clan $\GA$
in $(CO(X_0),\delta_{(X,X_0)})$ such that $F,G\in\GA$. Then, by the axiom (CS4), there exists a point $z\in X$ such that $\GA=\GA_{z,X_0}$. Therefore, $\GA=\{A\in CO(X_0)\st z\in\cl_X(A)\}$. Clearly, the principal filters $\FF$ and $\GG$ of the Boolean algebra $CO(X_0)$, which are generated, respectively, by $F$ and $G$, are contained in $\GA$. Since $\GA$ is a grill, Lemma \ref{grilllemma} implies that there exist ultrafilters $u$ and $v$ of the Boolean algebra $CO(X_0)$ such that $u\cup v\sbe \GA$, $\FF\sbe u$ and $\GG\sbe v$. Using the axiom (CS2), we obtain that there exist points $x,y\in X$ such that $\bigcap u=\{x\}$ and  $\bigcap v=\{y\}$. Then, clearly, $u=u_x$ and $v=u_y$. So, we get that $z\in\cl_X(A)$, for every $A\in u_x\cup u_y$. Thus, $xRy$. Obviously, $x\in F$ and $y\in G$. Therefore, $FC_RG$.
Conversely, let $F,G\in CO(X_0)$ and let $FC_RG$. Then $(\ex x\in F)(\ex y\in G)(xRy)$. Hence $F\in u_x$ and $G\in u_y$. Thus, $\cl_X(F)\cap\cl_X(G)\nes$. So, (\ref{xrycont}) is proved.
  This shows that the triple $(X,X_0,R)$ satisfies the axiom (PCS4). We obtain, as well, that $C_R=\d_{(X,X_0)}$ and thus, using  the axiom (CS4), we get that the axiom (PCS5) is satisfied as well.

Let us prove that the relation $R$ is a closed relation on the space $X_0$. Indeed, let $x,y\in X_0$ and $(x,y)\nin R$. Then, by (\ref{xrycont}),  there exist $F\in u_x$ and $G\in u_y$ such that $\cl_X(F)\cap\cl_X(G)=\ems$. This obviously implies that $(x,y)\in F\times G\sbe X\times X\stm R$. So, $R$ is a closed relation. Hence,
the triple $(X,X_0,R)$ satisfies the axiom (PCS2). Since the axioms (PCS1) and (PCS3) are obviously satisfied, we get that $(X,X_0,R)$ is a 2-precontact space.

 Let $(X,X_0,R\ap)$ be a 2-precontact space and $R\ap$ be a reflexive and symmetric relation on the set $X_0$. Then, by Proposition \ref{prop1}(b), $C_{R\ap}=(C_{R\ap})^\#$. Thus, using (\ref{clfg}), we get that for
every $F,G\in CO(X_0)$, $cl_X(F)\cap cl_X(G)\neq\emptyset \mbox{ iff }F(C_{R\ap}) G$. Now, (\ref{clxfclxgnescr}) implies that for
every $F,G\in CO(X_0)$, $F(C_R)G \iff F(C_{R\ap}) G$. Hence, $C_R\equiv  C_{R\ap}$. Set $B=CO(X_0)$. Then we get that $\unlb=(B,C_R)=(B,C_{R\ap})$. Using (\ref{corcx0}), the proof of Theorem \ref{mth}(b) and its notation, we obtain that: 1) $\Psi_2([(X,X_0,R\ap)])=[(B,C_R)]=[\unlb]$, 2) the map $f:(X,X_0,R\ap)\lra \unlx(B,C_R)$ is a PCS-isomorphism and, for every $x\in X_0$, $f(x)=u_x$. Let $\wh{R}=R_{(B,C_R)}$. Now, we get that for every $x,y\in X_0$, $xR\ap y\iff f(x)\wh{R} f(y)\iff u_x\wh{R}u_y\iff u_x\times u_y\sbe C_R\iff [(\fa F\in u_x)(\fa G\in u_y)(F(C_R)G)]\iff [(\fa F\in u_x)(\fa G\in u_y)(cl_X(F)\cap cl_X(G)\nes)]\iff xRy$. Therefore, $R\equiv R\ap$.
\sqs

\begin{cor}\label{cormthc}
If $\unlx=(X,X_0)$ is a 2-contact space then $X$ is a compact space.
\end{cor}

\doc By Lemma \ref{conprecon}, there exists a 2-precontact space $(X,X_0,R)$. Then, by Corollary \ref{cormth}, $X$ is a compact space.
\sqs

\begin{pro}\label{preconcon}
If $(X,X_0,R)$ is a reflexive and symmetric 2-precontact space then $(X,X_0)$ is a 2-contact space.
\end{pro}

\doc Clearly, the conditions (CS1)-(CS3) are fulfilled. Let $$\GA\in Clans(CO(X_0),\d_{(X,X_0)}).$$ By Proposition \ref{prop1}(b), we have that $C_R=(C_R)^\#$. Now (\ref{clfg}) implies that $C_R=\d_{(X,X_0)}$. Hence, $\GA\in Clans(CO(X_0),C_R)$ and, by (PCS5), there exists $x\in X$ such that $\GA=\GA_{x,X_0}$. So, the condition (CS4) is satisfied. Therefore $(X,X_0)$ is a 2-contact space.
\sqs

\begin{pro}\label{can2consp}
Let $\unlb=(B,C)$ be a contact algebra. Then the canonical 2-contact space $\unlxc(\unlb)$ is indeed a 2-contact space.
\end{pro}

\doc By Proposition \ref{proposition2}, the canonical 2-precontact space $\unlx(\unlb)=(X,X_0,R)$ is a 2-precontact space. Also, Lemma \ref{lmpc}(d)(e) implies that
$(X,X_0,R)$ is a reflexive and symmetric 2-precontact space. Now, using Proposition \ref{preconcon}, we get that $(X,X_0)$ is a 2-contact space. Since $\unlxc(\unlb)=(X,X_0)$, our assertion is proved.
\sqs

\begin{theorem}\label{mthc}\textsc{(New representation theorem for contact algebras.)}
\begin{description}
\item[{\rm (a)}]\quad  Let $(B,C)$ be a contact algebra and let
$(X,X_0)$ be the canonical 2-contact space of $(B,C)$ (see
Definition \ref{contactspace}(c)). Then the function $g_B:B\lra 2^X$, defined
in (\ref{fg}), is a CA-isomorphism from the algebra $(B,C)$ onto
the canonical contact algebra $(RC(X,X_0),C_{(X,X_0)})$ of
$(X,X_0)$. The sets $RC(X,X_0)$ and $RC(X)$ coincide iff the
contact algebra $(B,C)$ is complete.  The contact algebra $(B,C)$
is connected iff the 2-contact space $(X,X_0)$ is connected.

\item[{\rm (b)}]\quad There exists a bijective correspondence  between the class of
all, up to CA-isomor\-ph\-ism, (connected) contact algebras and the class of
all, up to CS-isomor\-phism, (connected) 2-contact spaces; namely,
for every CA $\unlb$, the CA-isomor\-phism class $[\unlb]$ of $\unlb$ corresponds to the CS-isomor\-phism class $[\unlx^c(\unlb)]$ of the canonical 2-contact space $\unlx^c(\unlb)$ of $\unlb$, and for every 2-contact space $(X,X_0)$, the CS-isomorphism class $[(X,X_0)]$ of $(X,X_0)$ corresponds to the CA-isomor\-phism class $[\unlb^c(X,X_0)]$ of the canonical contact algebra $\unlb^c(X,X_0)$ of $(X,X_0)$.
\end{description}
\end{theorem}

\doc (a) By Lemma \ref{conprecon}, there exists a unique reflexive and symmetric binary relation $R$ on $X_0$ such that
$(X,X_0,R)$ is a 2-precontact space. Since $(B,C)$ is a contact algebra, we have that $C=C^\#$. Thus, using (\ref{gisoc}), we get that
$g_B:(B,C)\lra (RC(X,X_0),C_{(X,X_0)})$ is a CA-isomorphism. The remaining assertions follow from Theorem \ref{mth}(a).

\medskip

\noindent(b) In this part of our proof, we will use the notation from the proof of Theorem \ref{mth}(b).

Let $\CC\AA$ be the class of all, up to CA-isomorphism, contact algebras. Let $\CC\SS$ be the class of all, up to CS-isomorphism, 2-contact spaces. Let
$2\CC\SS$ be the class of all 2-contact spaces. Let $2\PCS$ be the class of all 2-precontact spaces $(X,X_0,R)$ for which $R$ is a reflexive and symmetric relation.
Let $2\PP\SS$ be the class of all, up to  PCS-isomorphism, 2-precontact spaces $(X,X_0,R)$ for which $R$ is a reflexive and symmetric relation.
Then, using Lemma \ref{conprecon} and Proposition \ref{preconcon},
we get that the correspondence $$\p:2\CC\SS\lra 2\PCS,\ \ (X,X_0)\mapsto (X,X_0,R),$$ where the relation $R$ is defined by the formula (\ref{xrycont}), is a bijection. It is clear that then the correspondence $$\psi:\CC\SS\lra 2\PP\SS,\ \ [(X,X_0)]\lra [\p(X,X_0)],$$  is a bijection as well.
By Theorem \ref{mth}(a), $\unlb=(B,C)$ is a contact algebra iff the 2-precontact space $\unlx(\unlb)=(X,X_0,R)$ is reflexive and symmetric, i.e., iff $R$ is a reflexive and symmetric relation. Thus,  if $\Phi_2\ap$ is the restriction of the correspondence $\Phi_2$ to the subclass $\CC\AA$ of the class $\PCA$, then $$\Phi_2\ap:\CC\AA\lra 2\PP\SS$$ is a bijection. Therefore, the map $$\Phi_2^c=\psi\inv\circ\Phi_2\ap:\CC\AA\lra \CC\SS, \ \ [\unlb]\mapsto [\unlxc(\unlb)],$$ is a bijection.
The  assertion about connected contact algebras follows now from Theorem \ref{mth}(b).
\sqs

We are now going to obtain  an assertion from \cite{DV1} (namely, \cite[Theorem 5.1(ii)(for CAs)]{DV1})  as a corollary of Theorem \ref{mthc}. This assertion concerns the class of C-semiregu\-lar spaces introduced in \cite{DV1} (see Definition \ref{csemi} below). We will also obtain some new facts about this class of spaces. We start with recalling and proving some preliminary assertions. Then we obtain a new theorem about the structure of C-semiregular spaces (see Theorem \ref{csemith} below) and using it, we derive   \cite[Theorem 5.1(ii)(for CAs)]{DV1} from Theorem \ref{mthc} (see Corollary \ref{comcon} below).

\begin{lm}\label{dv141ii42i}{\rm (\cite{DV1})} Let $X$ be a topological space. Then:

\smallskip

\noindent(a) for every $x\in X$, $\s_x$ is a clan of the contact algebra $(RC(X),C_X)$;

\smallskip

 \noindent(b) if $X$ is semiregular, then $X$ is a $T_0$-space iff for every $x,y\in X$, $x\neq y$ implies that $\s_x\neq \s_y$ (see (\ref{sg}) for
$\sigma_x$).
 \end{lm}

 \begin{lm}\label{devrieslemma}{\rm (\cite{deV})}
 Let $(X,\TT)$ be a topological space and $U,V\in\TT$. Then $$\int(\cl(U)\cap \cl(V))=\int(\cl(U\cap V)).$$
 \end{lm}

 \begin{cor}\label{devrcor}
 Let $(X,\TT)$ be a topological space and $U,V\in\TT$. Then $$\cl(\int(\cl(U)\cap \cl(V)))=\cl(U\cap V).$$
 \end{cor}

 \doc Since $\cl(U\cap V)$ is a regular closed set, we get, using Lemma \ref{devrieslemma}, that $\cl(\int(\cl(U)\cap \cl(V)))=\cl(\int(\cl(U\cap V)))=\cl(U\cap V)$.
 \sqs

 \begin{defi}\label{upoint}
 \rm
 Let $(X,\TT)$ be a topological space and $x\in X$. The point $x$ is said to be an {\em u-point}\/ if for every $U,V\in\TT$, $x\in\cl(U)\cap \cl(V)$ implies that $x\in\cl(U\cap V)$.
 \end{defi}

\begin{pro}\label{remupoint}
 (a) A topological space $(X,\TT)$ is  extremally disconnected iff every of its points  is an u-point;

 \smallskip

 \noindent(b)   If $X$ is
 a dense subspace of a space $Y$ and $x\in X$, then $x$ is an u-point of $Y$ iff $x$ is an u-point of $X$.
\end{pro}

\doc (a) Let $X$ be extremally disconnected. Then, for every $U,V\in\TT$, we have, by Lemma \ref{devrieslemma}, that $\cl(U)\cap \cl(V)=\int(\cl(U)\cap \cl(V))=\int(\cl(U\cap V))=\cl(U\cap V)$. Hence, every  point of $X$  is an u-point.

Conversely, let every  point of $X$  be an u-point. Let $U\in\TT$. Suppose that $\cl (U)\nin\TT$. Then there exists $x\in\cl (U)\stm \int (\cl (U))$. Hence $x\in X\stm\int (\cl (U))=\cl (X\stm\cl (U))$. Set $V=X\stm\cl (U)$. Then $V\in\TT$ and $x\in\cl (U)\cap\cl (V)$. Since $x$ is an u-point, we get that $x\in\cl (U\cap V)=\ems$, a contradiction. Hence, $\cl(U)\in\TT$. So, $X$ is extremally disconnected.

\smallskip

 \noindent(b)  Let $x$ be an u-point of $X$. Let $U,V$ be open subsets of $Y$ and $x\in\cl_Y(U)\cap\cl_Y(V)$. Set $U\ap=U\cap X$ and $V\ap=V\cap X$. Then $x\in X\cap\cl_Y(U\ap)\cap\cl_Y(V\ap)=\cl_X(U\ap)\cap\cl_X(V\ap)$. Hence $x\in\cl_X(U\ap\cap V\ap)=\cl_X(X\cap U\cap V)\sbe \cl_Y(U\cap V)$. So, $x$ is an u-point of $Y$.

 Conversely, let $x$ be an u-point of $Y$. Let $U,V$ be open subsets of $X$ and $x\in\cl_X(U)\cap\cl_X(V)$. There exist open subsets $U\ap$ and $V\ap$ of $Y$ such that $U=U\ap\cap X$ and $V=V\ap\cap X$. Then $x\in\cl_Y(U\ap)\cap\cl_Y(V\ap)$. Hence $x\in X\cap \cl_Y(U\ap\cap V\ap)=X\cap \cl_Y(X\cap U\ap\cap V\ap)=X\cap\cl_Y(U\cap V)=\cl_X(U\cap V)$. Therefore, $x$ is an u-point of $X$.
\sqs

 \begin{pro}\label{proupoint}
 Let $(X,\TT)$ be a topological space and $x\in X$. Then $x$ is an u-point iff $\s_x$ is an ultrafilter of the Boolean algebra $RC(X)$ (see (\ref{sg}) for
$\sigma_x$).
 \end{pro}

 \doc Since, by Lemma \ref{dv141ii42i}(a), $\s_x$ is a grill of $RC(X)$, we have that $\s_x$ is an ultrafilter of $RC(X)$ iff $(\fa F,G\in \s_x)(F.G\in \s_x)$. Hence, using Corollary \ref{devrcor}, we get that:
 ($\s_x$ is an ultrafilter of $RC(X)$) $\iff$
 ($(\fa F,G\in RC(X))[(x\in F\cap G)\rightarrow(x\in F.G)]$) $\iff$
  ($(\fa F,G\in RC(X))[(x\in F\cap G)\rightarrow(x\in \cl(\int(F\cap G)))]$) $\iff$
   ($(\fa U,V\in\TT)[(x\in \cl(U)\cap \cl(V))\rightarrow(x\in \cl(\int(\cl(U)\cap \cl(V))))]$) $\iff$
    ($(\fa U,V\in\TT)[(x\in \cl(U)\cap \cl(V))\rightarrow(x\in \cl(U\cap V))]$) $\iff$
    ( $x$ is an u-point).
 \sqs

\begin{pro}\label{proupointopen}
 Let $X$ and $Y$ be  topological spaces and $x\in X$. If $x$ is an u-point and $f:X\lra Y$ is an open map, then $f(x)$ is an u-point.
 \end{pro}

 \doc Let $U$ and $V$ be  open subsets of $Y$ and $f(x)\in\cl_Y(U)\cap\cl_Y(V)$. We will show that $f(x)\in\cl_Y(U\cap V)$. We have that $x\in f\inv(\cl_Y(U))$ and
 $x\in f\inv(\cl_Y(V))$. Using \cite[Exercise 1.4.C]{E}, we get that $x\in \cl_X(f\inv(U))\cap\cl_X(f\inv(V))$. Since  $x$ is an u-point, we get, using again \cite[Exercise 1.4.C]{E}, that $x\in\cl_X(f\inv(U)\cap f\inv(V))=\cl_X(f\inv(U\cap V))=f\inv(\cl_Y(U\cap V))$. Thus, $f(x)\in\cl_Y(U\cap V)$. Hence $f(x)$ is an u-point.
 \sqs

Note that Proposition \ref{remupoint}(a) and Proposition \ref{proupointopen} imply the well-known fact that extremal disconnectedness is an invariant of open mappings (see \cite[Exercise 6.2.H(b)]{E}).

\begin{defi}\label{csemi}{\rm (\cite{DV1})}
\rm
A semiregular $T_0$-space $(X,\tau)$ is said to be \/ {\em C-semiregular}\/ if
for every clan $\Gamma$ in $(RC(X),C_X)$ there exists a point
$x\in X$ such that $\Gamma=\sigma_x$ (see (\ref{sg}) for
$\sigma_x$).
\end{defi}

The next assertion was stated in \cite{DV1} but it was left without proof there. For completeness, we will prove it here.

\begin{pro}\label{csemicomp}{\rm (\cite[Fact 4.1]{DV1})}
Every C-semiregular space $X$ is a compact space.
\end{pro}

\doc Let $\FF=\{F_a\st \a\in A\}$ be a centered (= with finite intersection property) family of closed subsets of $X$.
Since $X$ is semiregular, for every $\a\in A$ there exists a subfamily $\RR_\a$ of $RC(X)$ such that $F_\a=\bigcap\RR_\a$. Let $\RR=\bigcup\{\RR_\a\st \a\in A\}$. Then $\bigcap\FF=\bigcap\RR$ and $\RR$ is a centered family. Thus there exists an ultrafilter $u$ of the Boolean algebra $RC(X)$ containing $\RR$. Since $u$ is a clan in the contact algebra $(RC(X),C_X)$, there exists $x\in X$ such that $u=\s_x$. Therefore $x\in\bigcap u\sbe\bigcap \RR=\bigcap\FF$. Hence, $X$ is compact.
\sqs

\begin{theorem}\label{csemith}
 For every C-semiregular space $(X,\TT)$, the set $$X_0=\{x\in X\st x \mbox{ is an u-point of }X\}$$ endowed with its subspace topology is a dense extremally
 disconnected compact Hausdorff subspace of $(X,\TT)$ and is the unique dense extremally
 disconnected compact Hausdorff subspace of $(X,\TT)$.
 \end{theorem}

 \doc Set $(B,C)=(RC(X),C_X)$. Then $(B,C)$ is a complete Boolean algebra. Let $(\wh{X},\wh{X}_0)$ be the canonical 2-contact space of $(B,C)$ (see Definition \ref{contactspace}(c)). Since $B$ is a complete Boolean algebra, its Stone  space $\wh{X}_0$ is extremally disconnected and $RC(\wh{X},\wh{X}_0)=RC(\wh{X})$.
 Using Lemma \ref{dv141ii42i}(b) and the fact that $X$ is C-semiregular, we get that the map
 $$f:X\lra\wh{X},\ \ x\mapsto\s_x,$$
 is a bijection. In the notation of Definition
\ref{canspa}, we have that for every $F\in RC(X)$, $f(F)=\{f(x)\st x\in F\}=\{\s_x\st x\in F\}=\{\s_x\st F\in\s_x\}=\{\GA\in\wh{X}\st F\in\GA\}=g_B(F)$. Since $f$ is a bijection, we obtain now that for every $F\in RC(X)$,
$f\inv(g_B(F))=F$. Using the fact that $RC(X)$ is a closed base of $X$ and $\{g_B(F)\st F\in RC(X)\}$ is a closed base of $\wh{X}$, we get that $f$ is a homeomorphism. Since $\wh{X}_0=Ult(RC(X))$, we get that $\wh{X}_0=\{\s_x\st\s_x$ is an ultrafilter of $RC(X)\}$. Thus, Proposition \ref{proupoint} implies that $f\inv(\wh{X}_0)=X_0$. Hence, $f(X_0)= \wh{X}_0$. Therefore, $(X,X_0)$ is a 2-contact space and $f:(X,X_0)\lra(\wh{X},\wh{X}_0)$ is a CS-isomorphism. From this we get, in particular, that $X_0$ is  a
dense extremally disconnected compact Hausdorff subspace of $(X,\TT)$.

Let now $X_0\ap$ be a dense extremally disconnected compact Hausdorff subspace of $(X,\TT)$.
We will show that $X_0\ap=X_0$. Since $X_0\ap$ is extremally disconnected, Proposition \ref{remupoint}(a) implies that every point of $X_0\ap$ is an u-point of $X_0\ap$. Using the fact that $X_0\ap$ is dense in $X$, we get, by Proposition \ref{remupoint}(b), that every point of $X_0\ap$ is an u-point of $X$. Hence $X_0\ap\sbe X_0$.
 Obviously,  $X_0\ap$ is a dense subspace of $X_0$. Since $X_0$ is Hausdorff and $X_0\ap$ is compact, we get that $X_0\ap=X_0$.
 \sqs

 \begin{cor}\label{novocor}
 A compact Hausdorff space $X$ is C-semiregular iff it is extremally disconnected.
 \end{cor}

 \doc Let $X$ be extremally disconnected.. Then $RC(X)=CO(X)$. Also, the sets $Clans(RC(X),C_X)$ and $Ult(CO(X))$ coincide. Indeed, let $\GA\in Clans(RC(X),C_X)$. Then $\GA\sbe CO(X)$. Let $F,G\in\GA$. Since $F=(F\cap G)\cup(F\stm G)$ and $F\cap G,F\stm G\in CO(X)=RC(X)$, we have that $F\cap G\in\GA$ or $F\stm G\in \GA$. Clearly, $(F\stm G)(-C_X)G$. Hence $F\stm G\nin\GA$. Therefore, $F\cap G\in\GA$. Thus, $\GA$ is an ultrafilter of the Boolean algebra $CO(X)$. Then $|\bigcap\GA|=1$. Let $\{x\}=\bigcap\GA$. Then, obviously, $\GA=\s_x$. Hence, $X$ is a C-semiregular space.

 Conversely, if $X$ is C-semiregular then, by Theorem \ref{csemith}, $X$ contains a dense extremally disconnected compact Hausdorff space $Y$. Since $X$ is a Hausdorff space, we get that $X\equiv Y$. Therefore, $X$ is extremally disconnected.
 \sqs

\begin{cor}\label{novocor2}
 If $X$ is  C-semiregular  and $X_0=\{x\in X\st x \mbox{ is an u-point of }X\}$ then the pair $(X,X_0)$ is a 2-contact space and $X_0$ is a dense  extremally
 disconnected compact Hausdorff subspace of $X$; moreover, $X_0$ is the unique dense  extremally
 disconnected compact Hausdorff subspace of $X$.
 \end{cor}

 \doc By Theorem \ref{csemith}, $X_0$ is a dense extremally
 disconnected compact Hausdorff subspace of $X$ and hence, the conditions (CS1) and (CS2) are fulfilled. Also, we obtain that $RC(X_0)=CO(X_0)$ and thus $RC(X,X_0)=RC(X)$; moreover, the map $$e:(CO(X_0),\d_{(X,X_0)})\lra (RC(X),C_X), \ \ F\mapsto\cl_X(F),$$ is a CA-isomorphism. Since $X$ is semiregular, we get that the condition (CS3) is fulfilled. Let $\GA\in Clans(CO(X_0),\d_{(X,X_0)})$. Then $e(\GA)\in Clans(RC(X),C_X)$ and, therefore, there exists $x\in X$ such that $e(\GA)=\s_x$. Since, obviously, $\s_x=\GA_{x,X_0}$, we get that the condition (CS4) is also fulfilled. Hence $(X,X_0)$ is a 2-contact space.

 The uniqueness assertion follows from Theorem \ref{csemith}.
 \sqs

\begin{lm}\label{csem2cont}
If $(X,X_0)$ is a 2-contact space and $X_0$ is extremally disconnected, then $X$ is C-semiregular.
\end{lm}

\doc By (CS1), $X$ is a $T_0$-space. Using Lemma \ref{lmrcx0}(a) and (CS3), we get that $X$ is a semiregular space and $RC(X,X_0)=RC(X)$. Since $CO(X_0)=RC(X_0)$, Lemma \ref{isombool} implies that the map
$$e:(CO(X_0),\d_{(X,X_0)})\lra (RC(X),C_X),\ \ F\mapsto\cl_X(F),$$
is a CA-isomorphism. Let $\GA\in Clans(RC(X),C_X)$ and $\GA\ap=e\inv(\GA)$. Then $\GA\ap\in Clans(CO(X_0),\d_{(X,X_0)})$. Thus, by (CS4), there exists $x\in X$ such that $\GA\ap=\GA_{x,X_0}$. Then we get that $\GA=e(\GA\ap)=e(\GA_{x,X_0})=\s_x$. Therefore, $X$ is a C-semiregular space.
\sqs

\begin{cor}\label{comcon}{\rm (\cite{DV1})}
There exists a bijective correspondence  between the class of all,
up to CA-isomor\-ph\-ism, (connected) complete  contact algebras
and the class of all, up to homeomorphism, (connected)
C-semiregular spaces.
\end{cor}

\doc Let $\CC\CC\AA$ be the class of all, up to CA-isomorphism, complete contact algebras. Let $\CC\CC\SS$ be the class of all, up to CS-isomorphism, 2-contact spaces $(X,X_0)$ such that $X_0$ is extremally disconnected.  Using Theorem \ref{mthc}(b) and the notation of its proof, we get that the map
$$\Phi_2^c:\CC\AA\lra \CC\SS, \ \ [(B,C)]\mapsto [(X,X_0)],$$
where $(X,X_0)$ is the canonical 2-precontact space of the contact algebra $(B,C)$ (see Definition \ref{contactspace}(c)), is a bijection. Then Theorem \ref{mthc}(a) implies that the restriction $\Phi_c$ of the correspondence $\Phi_2^c$ to the class $\CC\CC\AA$ is a bijection between the later class and the class $\CC\CC\SS$, i.e.,
$$\Phi_c:\CC\CC\AA\lra \CC\CC\SS$$
is a bijection.
Let $\CC\SS\RR$ be the class of all, up to homeomorphism, C-semiregular spaces. Then the map
$$\a:\CC\CC\SS\lra \CC\SS\RR,\ \ [(X,X_0)]\mapsto [X],$$
is a bijection. Indeed, Lemma \ref{csem2cont} implies that the correspondence $\a$ is well-defined; the fact that $\a$ is a bijection follows from Theorem \ref{csemith}.
Now we get that the composition $$\Phi_c\ap=\a\circ\Phi_c:\CC\CC\AA\lra \CC\SS\RR$$ is a bijection.
The  assertion about connected complete contact algebras follows now from Theorem \ref{mth}(b).
\sqs

\section{A connected version of the Stone Duality Theorem}

The new representation theorems, presented in the previous section, permit us to obtain as particular cases the Stone Representation Theorem \cite{ST} and a new connected version of it.
Let us start with the Stone Representation Theorem \cite{ST}.

\begin{pro}\label{stonerepth}
The Stone Representation Theorem follows from Theorem \ref{mthc}; namely, we obtain it equipping each Boolean algebra $B$ with the smallest contact relation $\rho_s^B$ on $B$ (see Example \ref{extrcr} for $\rho_s^B$).
\end{pro}

\doc Let $B$ be a Boolean algebra. Then, by \cite[Corollary 3.3]{DV1}, $Ult(B,\rho_s)=Clans(B,\rho_s)$. Hence, the canonical 2-contact space of the contact algebra $\unlb=(B,\rho_s)$ is the pair $\unlx^c(\unlb)=(X,X_0)$, where $X=X_0$ is the Stone space of $B$. Also, $RC(X,X_0)=CO(X_0)$ and $C_{(X,X_0)}=\rho_s^{B\ap}$, where $B\ap=CO(X_0)$. So, Theorem \ref{mthc}(a) reduces to the Stone Representation Theorem. Note that $RC(X_0)=CO(X_0)$ iff $X_0$ is extremally disconnected. Note as well that $(B,\rho_s)$ is connected iff $B=\{0,1\}$. Further, Theorem \ref{mthc}(b) reduced to the Stone Theorem that there exists a bijective correspondence  between
 the class of all, up to Boolean isomor\-ph\-ism, Boolean algebras and the class of
all, up to homeomorphism, Stone spaces. Finally, Corollary \ref{comcon} and Corollary \ref{novocor} imply that there exists a bijective correspondence  between
 the class of all, up to Boolean isomor\-ph\-ism, complete Boolean algebras and
 the class of all, up to homeomorphism, compact Hausdorff extremally disconnected spaces.
\sqs

Now, equipping each Boolean algebra $B$ with the largest contact relation $\rho_l^B$ on $B$ (see Example \ref{extrcr} for $\rho_l^B$), we will obtain a connected version of the Stone Representation Theorem. Further on we will extend it to a connected version of the Stone Duality Theorem.

\begin{defi}\label{2grillsp}\textsc{(Stone 2-spaces.)}
\rm

\smallskip

\noindent(a) A topological pair  $(X,X_0)$ is  called a\/ {\em Stone
2-space} (abbreviated as S2S) if it satisfies conditions (CS1)-(CS3) of Definition \ref{contactspace} and the following condition:

\begin{description}
\item[$(S2S4)$]\quad If $\Gamma\in
Grills(CO(X_0))$  then there exists a point
$x\in X$ such that $\Gamma =\Gamma_{x,X_0}$ (see (\ref{sg}) for
$\Gamma_{x,X_0}$).
\end{description}

\smallskip

\noindent(b) Let   $(X,X_{0})$ and $(\widehat{X},\widehat{X}_{0})$ be two
Stone 2-spaces. We say that  $(X,X_0)$ and
$(\widehat{X},\widehat{X}_0)$ are\/ {\em S2S-isomorphic} (or,
simply,\/ {\em isomorphic}) if there exists a homeomorphism
$f:X\longrightarrow\widehat{X}$ such that
  $f(X_{0})=\widehat{X}_0$.
\end{defi}

\begin{pro}\label{constonereppro}
Let $(X,X_0)$ be a Stone 2-space and $B=RC(X,X_0)$. Then:

\smallskip

\noindent(a) $C_{(X,X_0)}=\rho_l^B$  (see Example \ref{extrcr} for $\rho_l^B$), and

\smallskip

\noindent(b)  $(X,X_0)$ is a 2-contact space.
\end{pro}

\doc (a)  Clearly, $C_{(X,X_0)}\sbe\rho_l^B$. Recall that for any $F,G\in CO(X_0)$, $$ \cl_X(F) C_{(X,X_0)} \cl_X(G)\iff \cl_X(F)\cap\cl_X(G)\nes.$$
 We will show that $\rho_l^B\sbe C_{(X,X_0)}$. Let $F,G\in CO(X_0)$ and $F\nes,\ G\nes$ (i.e., $\cl_X(F),\cl_X(G)\in B$ and $\cl_X(F) \rho_l^B\cl_X(G)$). Obviously, there exist $u,v\in Ult(CO(X_0))$ such that $F\in u$ and $G\in v$. Then $\GA=u\cup v\in Grills(CO(X_0))$ (see, e.g., \cite[Corollary 3.1]{DV1}). Hence, by (S2S4), there exists a point
$x\in X$ such that $\Gamma =\Gamma_{x,X_0}$. Then $x\in \cl_X(F)\cap\cl_X(G)$ and thus $\cl_X(F) C_{(X,X_0)} \cl_X(G)$. So, $\rho_l^B= C_{(X,X_0)}$.

\smallskip

\noindent(b) Since the map $e:(CO(X_0),\d_{(X,X_0)})\lra (RC(X,X_0),C_{(X,X_0)}), \ \ F\mapsto\cl_X(F),$ is a CA-isomorphism, we get, using (a) and \cite[Example 3.1]{DV1}, that $$Clans(CO(X_0),\d_{(X,X_0)})=Grills(CO(X_0)).$$ Thus, condition (CS4) follows from condition (S2S4). Hence  $(X,X_0)$ is a 2-contact space.
\sqs

\begin{cor}\label{conn2grsp}
Let $(X,X_0)$ be a Stone 2-space. Then $X$ is a compact connected $T_0$-space.
\end{cor}

\doc According to Proposition \ref{constonereppro}(b), $(X,X_0)$ is a 2-contact space. Hence, by Corollary \ref{cormthc}, $X$ is a compact space. Since the canonical CA $(RC(X,X_0),C_{(X,X_0)})$ of the 2-contact space $(X,X_0)$ is connected (by Proposition \ref{constonereppro}(a)), Theorem \ref{mthc}(a) implies that the space $X$ is connected.
\sqs

\begin{defi}\label{2grillcansp}\textsc{(The canonical Stone 2-space of a Boolean algebra and the canonical Boolean algebra of a Stone 2-space.)}
\rm

\smallskip

\noindent(a) Let $B$ be a Boolean algebra. Then the canonical 2-contact space of the contact algebra $(B,\rho_l^B)$ is said to be the {\em canonical Stone 2-space of the Boolean algebra} $B$ and is denoted by $\unlx^s(B)$.

\smallskip

\noindent(b) Let $(X,X_0)$ be a Stone 2-space. Then the Boolean algebra $RC(X,X_0)$ is said to be the {\em canonical  Boolean algebra of} $(X,X_0)$.
\end{defi}

\begin{pro}\label{can2grsp}
The canonical Stone 2-space of a Boolean algebra $B$ is indeed a Stone 2-space.
\end{pro}

\doc Let $\unlx^c(B,\rho_l^B)=(X,X_0)$. Then, by Theorem \ref{mthc}(a), the contact algebras $(B,\rho_l^B)$ and $(RC(X,X_0),C_{(X,X_0)})$ are CA-isomorphic. Since the contact algebras $(RC(X,X_0),C_{(X,X_0)})$ and $(CO(X_0),\d_{(X,X_0)})$ are CA-isomorphic, we get that the contact algebras $(B,\rho_l^B)$ and $(CO(X_0),\d_{(X,X_0)})$ are CA-isomorphic. Then, by \cite[Example 3.1]{DV1}, $Clans(CO(X_0),\d_{(X,X_0)})=Grills(CO(X_0))$. According to Proposition \ref{can2consp}, $(X,X_0)$ is a 2-contact space. Now we obtain that $(X,X_0)$ is a Stone 2-space.
\sqs

\begin{theorem}\label{constoneth}\textsc{(A connected version of the Stone Representation Theorem.)}
\begin{description}
\item[{\rm (a)}]\quad  Let $B$ be a Boolean algebra and let
$(X,X_0)$ be the canonical Stone 2-space of $B$ (see
Definition \ref{2grillcansp}). Then the function $g_B:B\lra 2^X$, defined
in (\ref{fg}), is a Boolean isomorphism from the Boolean algebra $B$ onto
the canonical Boolean algebra $RC(X,X_0)$ of
$(X,X_0)$. The sets $RC(X,X_0)$ and $RC(X)$ coincide iff the
Boolean algebra $B$ is complete.

\item[{\rm (b)}]\quad There exists a bijective correspondence  between the class of
all, up to Boolean isomor\-ph\-ism, Boolean algebras and the class of
all, up to S2S-isomorphism,  Stone 2-spaces; namely, for every Boolean algebra $B$, $[B]\mapsto [\unlx^s(B)]$, and for every Stone 2-space $(X,X_0)$, $[(X,X_0)]\mapsto [RC(X,X_0)]$.
\end{description}
\end{theorem}

\doc (a) Since $(X,X_0)$ is, by Definition \ref{2grillcansp}(a), the canonical 2-contact space of the contact algebra $(B,\rho_l^B)$, our assertion follows from Theorem \ref{mthc}(a).

\smallskip

\noindent(b)  It follows from Theorem \ref{mthc}(b), Definition \ref{2grillcansp}(a), Proposition \ref{can2grsp} and Proposition \ref{constonereppro}.
\sqs

\begin{defi}\label{grillsp}\textsc{(Extremally connected spaces.)}
\rm
 A semiregular $T_0$-space $X$ is  called an\/ {\em extremally connected
space} (abbreviated as ECS) if for every grill $\GA$ in $RC(X)$ there exists a point
$x\in X$ such that $\Gamma=\sigma_x$ (see (\ref{sg}) for
$\sigma_x$).
\end{defi}

\begin{pro}\label{grspCsemrs}
Let $X$ be an extremally connected space and $B=RC(X)$. Then:

\smallskip

\noindent(a) $C_X=\rho_l^B$ (see Lemma \ref{rcl} for $C_X$ and Example \ref{extrcr} for $\rho_l^B$), and

\smallskip

\noindent(b) $X$ is a C-semiregular space.
\end{pro}

\doc (a) Let $F,G\in B$ and $F\nes, \ G\nes$. Then, as in the proof of Proposition \ref{constonereppro}(a), we get that there exists a grill $\GA$ in $B$ such that
$F,G\in\GA$. Then, by Definition \ref{grillsp}, there exists $x\in X$ such that $\GA=\s_x$. Thus $x\in F\cap G$, i.e. $FC_X G$. Therefore, $\rho_l^B\sbe C_X$. Since, obviously, $C_X\sbe \rho_l^B$, we get that $C_X=\rho_l^B$.

\smallskip

\noindent(b) Since, by (a), $(B,C_X)=(B,\rho_l^B)$, \cite[Example 3.1]{DV1} implies that $Clans(B,C_X)=Grills(B)$.  Thus $X$ is a C-semiregular space.
\sqs

\begin{cor}\label{grsp2grsp}
Let $X$ be an extremally connected space and $$X_0=\{x\in X\st x \mbox{ is an u-point of }X\}.$$ Then $(X,X_0)$ is a Stone 2-space and $X_0$ is a dense  extremally disconnected compact Hausdorff subspace of $X$; moreover, $X_0$ is the unique  dense  extremally disconnected compact Hausdorff subspace of $X$.
\end{cor}

\doc According to Proposition \ref{grspCsemrs}(b), $X$ is a C-semiregular space. Then, by Corollary \ref{novocor2}, $(X,X_0)$ is a 2-contact space and $X_0$ is a dense extremally disconnected compact Hausdorff subspace of $X$. Let $\GA\in Grills (CO(X_0))$. We have that $CO(X_0)=RC(X_0)$ and thus, by Lemma \ref{isombool},
the Boolean algebra $CO(X_0)$ is isomorphic to the Boolean algebra $RC(X)$. Now, using \cite[Example 3.1]{DV1} and
Proposition \ref{grspCsemrs}(a), we get that $Grills(CO(X_0))=Clans(CO(X_0),\rho_l^{CO(X_0)})=Clans(CO(X_0),\d_{(X,X_0)})$. Hence, by condition (CS4), there exists $x\in X$ such that $\GA=\GA_{x,X_0}$. Therefore, condition (S2S4) is satisfied and we obtain
 that $(X,X_0)$ is a Stone 2-space.

 The uniqueness assertion follows from Theorem  \ref{csemith}.
\sqs

\begin{lm}\label{2grgr}
If $(X,X_0)$ is a Stone 2-space and $X_0$ is extremally disconnected, then $X$ is an extremally connected space.
\end{lm}

\doc Let $B=RC(X)$. Since $X_0$ is an extremally disconnected  space, we have that $RC(X_0)=CO(X_0)$. Then $RC(X,X_0)=B$. Now, Proposition \ref{constonereppro}(a) implies that $C_X=\rho_l^B$ and $(X,X_0)$ is a 2-contact space. Then, by Lemma \ref{csem2cont}, $X$ is a C-semiregular space. Since $(B,C_X)=(B,\rho_l^B)$, \cite[Example 3.1]{DV1} implies that $X$ is an extremally connected space.
\sqs

\begin{theorem}\label{constonethcomp}
There exists a bijective correspondence  between the class of
all, up to Boolean isomor\-ph\-ism, complete Boolean algebras and the class of
all, up to homeomorphism,  extremally connected spaces; namely, for every complete Boolean algebra $B$, $[B]\mapsto [X]$, where $X$ is the first component of $\unlx^s(B)$, and for every extremally connected space $X$, $[X]\mapsto [RC(X)]$.
\end{theorem}

\doc It follows from Theorem \ref{constoneth}(b), the Stone bijection between complete Boolean algebras and extremally disconnected Stone spaces, Lemma \ref{2grgr} and Corollary \ref{grsp2grsp}.
\sqs

Now we will obtain a connected version of the Stone Duality Theorem. We start with three simple assertions, the first of which is a slight generalization of \cite[Proposition 4.1(iv)]{DV1}.

\begin{lm}\label{pro41dv1}
Let $X$ be a topological space, $B$ be a subalgebra of the Boolean algebra $(RC(X),+,.,*,\ems,X)$ (defined in \ref{rcdef}), $\GA$ be a grill of $B$, $x\in X$ and $\GA\sbe\s_x^B$. Then $\nu_x^B\sbe \GA$ (see (\ref{sg1}) for $\nu_x^B$ and (\ref{sg}) for $\s_x^B$).
\end{lm}

\doc We can suppose that $\GA\neq\s_x^B$. Let $F\in\s_x^B\stm\GA$. Since $F+F^*=1\in\GA$, we get that  $F^*\in\GA$. Suppose that $F\in\nu_x^B$. Then $x\in\int(F)$ and thus $x\nin F^*$. Since $\GA\sbe\s_x^B$, we get a contradiction. Therefore, $F\nin\nu_x^B$. So, we obtain that $\nu_x^B\cap(\s_x^B\stm\GA)=\ems$. Since $\nu_x^B\sbe\s_x^B$, we get that $\nu_x^B\sbe\GA$.
\sqs

\begin{lm}\label{lmconst}
Let $X$ be a topological space, $B$ be a subalgebra of the Boolean algebra $(RC(X),+,.,*,\ems,X)$ (defined in \ref{rcdef}) and $B$ be a closed base for the space $X$. Then, for every $x\in X$, the family $\BB_x=\{\int(F)\st F\in\nu_x^B\}$ is a base for $X$ at the point $x$.
\end{lm}

\doc Since $B$ is a closed base for $X$, the family $\BB=\{X\stm F\st F\in B\}$ is an open base for $X$. Further, for any $F\in B$, we have that $\int(F^*)=X\stm F$. Thus $\BB=\{\int(F^*)\st F\in B\}=\{\int(F)\st F\in B\}$. Hence, for every $x\in X$, the family $\{U\in\BB\st x\in U\}$ is a base for $X$ at the point $x$. Clearly,
$\{U\in\BB\st x\in U\}=\{\int(F)\st F\in B, x\in\int(F)\}=\{\int(F)\st F\in \nu_x^B\}=\BB_x$. Therefore, for every $x\in X$, the family $\BB_x$ is a base for $X$ at the point $x$.
\sqs

\begin{lm}\label{lmgrrestr}
Let $A$ be a subalgebra of a Boolean algebra $B$ and $\GA\in Grills(B)$. Then $\GA\cap A\in Grills(A)$.
\end{lm}

\doc Clearly, if $u\in Ult(B)$ then $u\cap A\in Ult(A)$. Then, using \cite[Corollary 3.1]{DV1}, we obtain our assertion.
\sqs

\begin{defi}\label{defeconsp}
\rm
Let $(X,X_0)$ and $(X\ap,X_0\ap)$ be two Stone 2-spaces and $f:X\lra X\ap$ be a continuous map. Then $f$ is called a {\em 2-map}\/ if $f(X_0)\sbe X_0\ap$.

The category of all Stone 2-spaces and all 2-maps between them will be denoted by $\ECS$.

The category of all Boolean algebras and all Boolean homomorphisms between them will be denoted by $\Boo$.
\end{defi}

\begin{theorem}\label{constoneduality}
The categories $\Boo$ and $\ECS$ are dually equivalent.
\end{theorem}

\doc We will first define two
contravariant functors $$D^a:\Boo\lra\ECS\ \ \mbox{ and }\ \
D^t:\ECS\lra\Boo.$$

Let $(X,X_0)\in\card{\ECS}$. Define $$D^t(X,X_0)=RC(X,X_0),$$
i.e. $D^t(X,X_0)$ is the canonical Boolean algebra of the Stone 2-space $(X,X_0)$ (see Definition \ref{2grillcansp}(b)).
 Hence $D^t(X,X_0)\in\card{\Boo}$.

Let $f\in\ECS((X,X_0),(Y,Y_0))$. Define $D^t(f):D^t(Y,Y_0)\lra D^t(X,X_0)$ by
the formula
\begin{equation}\label{deftetafc}
D^t(f)(\cl_Y(G))=\cl_X(X_0\cap f\inv(G)), \ \  \fa G\in CO(Y_0).
\end{equation}
Set $\p_f=D^t(f)$. We will show that $\p_f$ is a Boolean homomorphism
between the Boolean algebras $RC(Y,Y_0)$ and $RC(X,X_0)$. Clearly, $\p_f(\ems)=\ems$ and $\p_f(Y)=X$.
Let $F,G\in CO(Y_0)$. Then $\p_f(\cl_Y(F)+\cl_Y(G))=\p_f(\cl_Y(F\cup G))=\cl_X(X_0\cap f\inv(F\cup G))=\cl_X((X_0\cap f\inv(F))\cup(X_0\cap f\inv(G)))=\p_f(\cl_Y(F))+
\p_f(\cl_Y(G))$. Also, using Lemma \ref{isombool}, we get that $\p_f((\cl_Y(F))^*)=\p_f(\cl_Y(Y_0\stm F)=\cl_X(X_0\cap f\inv(Y_0\stm F))=\cl_X(X_0\cap(f\inv(Y_0)\stm f\inv(F)))=\cl_X(X_0\stm(X_0\cap f\inv(F)))=(\cl_X(X_0\cap f\inv(F)))^*=(\p_f(\cl_Y(F))^*$. So, $\p_f$ is a Boolean homomorphism, i.e. $D^t(f)$ is well-defined.

Now we will show  that $D^t$ is a contravariant functor. Clearly, $D^t(id_{X,X_0})=id_{D^t(X,X_0)}$. Let $f\in\ECS((X,X_0),(Y,Y_0))$ and $g\in\ECS((Y,Y_0),(Z,Z_0))$.
Then, for every $F\in CO(Z_0)$, $D^t(g\circ f)(\cl_Z(F))=\cl_X(X_0\cap(g\circ f)\inv(F))=\cl_X(X_0\cap f\inv(g\inv(F)))$ and $(D^t(f)\circ D^t(g))(\cl_Z(F))=D^t(f)(\cl_Y(Y_0\cap g\inv(F)))=\cl_X(X_0\cap f\inv(Y_0\cap g\inv(F)))=\cl_X(X_0\cap f\inv(Y_0)\cap f\inv(g\inv(F)))=\cl_X(X_0\cap f\inv(g\inv(F)))=
D^t(g\circ f)(\cl_Z(F))$. So, $D^t$ is a contravariant functor.

For every Boolean algebra  $A$, set $$D^a(A)=(X,X_0),$$
where $(X,X_0)$ is the canonical Stone 2-space of the Boolean algebra $A$ (see Definition \ref{2grillcansp}(a)).
Then Proposition \ref{can2grsp} implies that
$D^a(A)\in\card{\ECS}$.

Let $\p\in\Boo(A,B)$. Let $D^a(A)=(X,X_0)$ and $D^a(B)=(Y,Y_0)$.
Then we define the map
$$D^a(\p):D^a(B)\lra D^a(A)$$ by the formula
\begin{equation}\label{deftetaphic}
D^a(\p)(\GA)=\p\inv(\GA), \ \  \fa \GA\in Y.
\end{equation}
Set $f_\p=D^a(\p)$. Since every grill of a Boolean algebra $B\ap$ is a union of ultrafilters of $B\ap$ and every union of ultrafilters of $B\ap$ is a grill of $B\ap$ (see, e.g., \cite[Corollary 3.1]{DV1}), and the inverse image of an ultrafilter by a Boolean homomorphism between two Boolean algebras is again an ultrafilter, we get that
$\fa \GA\in Y$, $f_\p(\GA)\in X$, i.e.
 $f_\p:Y\lra X$.

We will show that $f_\p$ is a continuous function. Let $a\in A$. Then $g_A(a)=\{\GA\in X\st a\in \GA\}$ is a basic closed subset of $X$ (see Definition \ref{2grillcansp}(a)). We will show that
\begin{equation}\label{contfphic}
f_\p\inv(g_A(a))=g_B(\p(a))(=\{\GA\ap\in Y\st\p(a)\in\GA\ap\}).
\end{equation}
Indeed, let $\GA\ap\in f_\p\inv(g_A(a))$. Then $f_\p(\GA\ap)\in g_A(a)$. Thus $a\in
\p\inv(\GA\ap)$, i.e. $\p(a)\in\GA\ap$. So, $\GA\ap\in g_B(\p(a))$. Hence $f_\p\inv(g_A(a))\sbe g_B(\p(a))$. Conversely,
let $\GA\ap\in g_B(\p(a))$, i.e. $\p(a)\in\GA\ap$. Then $a\in\p\inv(\GA\ap)=f_\p(\GA\ap)$. Hence $f_\p(\GA\ap)\in g_A(a)$. Then $\GA\ap\in f_\p\inv(g_A(a))$. So,
$f_\p\inv(g_A(a))\spe g_B(\p(a))$. Thus the equation (\ref{contfphic}) is verified and
we get that $f_\p$ is a continuous function.

Let us now show that $f_\p(Y_0)\sbe X_0$. Let $u\ap\in Y_0$. Then $u\ap\in Ult(B)$. Hence $f_\p(u\ap)=\p\inv(u\ap)\in Ult(A)=X_0$. Therefore, $f_\p(Y_0)\sbe X_0$.
So,
$$D^a(\p)\in\ECS(D^a(B),D^a(A)).$$

Clearly, for every Boolean algebra $B$, $D^a(id_B)=id_{D^a(B)}$. Let $\p\in\Boo(A,B)$ and $\psi\in\Boo(B,B\ap)$. Let $f_\p=D^a(\p)$, $f_\psi=D^a(\psi)$ and $D^a(B\ap)=(Z,Z_0)$. Then, for every $\GA\in Z$, we have that $D^a(\psi\circ\p)(\GA)=(\psi\circ\p)\inv(\GA)=\p\inv(\psi\inv(\GA))=f_\p(f_\psi(\GA))=(D^a(\p)\circ D^a(\psi))(\GA)$. We get
 that $D^a$ is a contravariant functor.

Let $(X,X_0)\in\card{\ECS}$. Then $D^t(X,X_0)=RC(X,X_0)$. Set $B=RC(X,X_0)$. Let $D^a(B)=(Y,Y_0)$. Then $Y=Grills(B)$ and $Y_0=Ult(B)$.
By \cite[Proposition 4.1(ii)]{DV1}, if $x\in X$ then $\s_x\in Clans(RC(X),C_X)$. Using Lemma \ref{lmgrrestr}, we get that, for every $x\in X$, $\s_x^B\in Clans(B,C_{(X,X_0)})$. According to Proposition \ref{constonereppro}(a), $C_{(X,X_0)}=\rho_l^B$  (see Example \ref{extrcr} for $\rho_l^B$). Hence, by \cite[Example 3.1]{DV1}, $Clans(B,C_{(X,X_0)})=Grills(B)$. Therefore, for every $x\in X$, $\s_x^B\in Grills(B)$.
So, the following map is well-defined:
\begin{equation}\label{txcnc}
t_{(X,X_0)}:(X,X_0)\lra D^a(D^t(X,X_0)), \ \ x\mapsto \s_x^{RC(X,X_0)},
\end{equation}
and we will  show that it is a homeomorphism.  We start by proving that $t_{(X,X_0)}$ is a surjection. Let $\GA\in Y$. Then $\GA\in Grills(B)$. Using Lemma \ref{isombool}, we get that $\GA\ap=r_{X_0,X}(\GA)\in Grills(CO(X_0))$. Hence, by (S2S4), there exists $x\in X$ such that $\GA\ap=\GA_{x,X_0}$. Since, by Lemma  \ref{isombool}, $\GA=e_{X_0,X}(\GA\ap)$, we get that $\GA=\s_x^B=t_{(X,X_0)}(x)$. So, $t_{(X,X_0)}$ is a surjection. For showing that $t_{(X,X_0)}$ is a injection, let $x,y\in X$ and $x\neq y$. Since $X$ is a $T_0$-space, there exists an open subset $U$ of $X$ such that $|U\cap\{x,y\}|=1$. We can suppose,
 without loss of generality, that $x\in U$ and $y\nin U$. Since $B$ is a closed base of $X$, there exists $F\in B$ such that $x\in X\stm F\sbe U$. Then $y\in F$ and $x\nin F$. Hence $F\in \s_y^B$ and $F\nin\s_x^B$, i.e. $t_{(X,X_0)}(x)\neq t_{(X,X_0)}(y)$. So, $t_{(X,X_0)}$ is a injection. Thus $t_{(X,X_0)}$ is a bijection.
We will now prove that $t_{(X,X_0)}$ is a continuous map. We have that the family $\{g_B(F)=\{\GA\in Y\st F\in\GA\}\st F\in B\}$ is a closed base of $Y$. Let $F\in B$. We will show that
\begin{equation}\label{txcontc}
t_{(X,X_0)}\inv(g_B(F))=F.
\end{equation}
Let $x\in F$. Set $t_{(X,X_0)}(x)=\GA$. Then $\GA=\s_x^B$. Since $F\in\GA$, we get that $\GA\in g_B(F)$. Thus $t_{(X,X_0)}(F)\sbe g_B(F)$, i.e. $F\sbe t_{(X,X_0)}\inv(g_B(F))$. Conversely, let $x\in t_{(X,X_0)}\inv(g_B(F))$. Set $\GA=t_{(X,X_0)}(x)$. Then $\GA\in g_B(F)$. Hence $F\in\GA$. Since $\GA=\s_x^B$, we get that $x\in F$. Hence $F\spe t_{(X,X_0)}\inv(g_B(F))$. So, $F= t_{(X,X_0)}\inv(g_B(F))$. This shows that $t_{(X,X_0)}$ is a continuous map. For showing that $t_{(X,X_0)}\inv$ is a continuous map, let $F\in B$. Using (\ref{txcontc}) and the fact that $t_{(X,X_0)}$ is a bijection, we get that $t_{(X,X_0)}(F)=g_B(F)$. Hence $(t_{(X,X_0)}\inv)\inv(F)=g_B(F)$. This shows that $t_{(X,X_0)}\inv$ is a continuous map. So, $t_{(X,X_0)}$ is a homeomorphism.

We will now  show that $t_{(X,X_0)}(X_0)=Y_0$. Let $x\in X_0$. Set $\GA=t_{(X,X_0)}(x)$. Then $\GA=\s_x^B$ and $r_{X_0,X}(\GA)=\{F\in CO(X_0)\st x\in F\}=u_x^{X_0}\in Ult(CO(X_0)$. Then, by Lemma \ref{isombool}, $\GA=e_{X_0,X}(u_x^{X_0})\in Ult(B)=Y_0$. Hence $t_{(X,X_0)}(X_0)\sbe Y_0$. Let now $\GA\in Y_0$. Then $\GA\in Ult(B)$ and thus $u=r_{X_0,X}(\GA)\in Ult(CO(X_0))$. Clearly, there exist $x\in X_0$ such that $u=u_x^{X_0}$. Then $\GA=e_{X_0,X}(u_x^{X_0})=\s_x^B=t_{(X,X_0)}(x)$. Therefore, $t_{(X,X_0)}(X_0)\spe Y_0$. We have proved that $t_{(X,X_0)}(X_0)= Y_0$. So, $t_{(X,X_0)}$ is a $\ECS$-isomorphism.

Let $B$ be a Boolean algebra and let us set $(X,X_0)=D^a(B)$. Then $D^t(X,X_0)=RC(X,X_0)$ and, using Theorem \ref{constoneth}(a), we get that the map $$g_B:B\lra RC(X,X_0), \ \  a\mapsto g_B(a)=\{\GA\in X\st a\in\GA\},$$
is a Boolean isomorphism.

We will now show that
%
$$t:Id_{\ECS}\lra   D^a\circ  D^t,$$
%
 defined by
%
 $t(X,X_0)=\tcx0, \ \ \fa (X,X_0)\in\card\ECS,$
%
is a natural isomorphism.

Let $f\in\ECS((X,X_0),(Y,Y_0))$ and $\fs= D^a(D^t(f))$. We have to show
that $\fs\circ \tcx0=\tcy0\circ f$. Set $\p_f=D^t(f)$, $A=RC(X,X_0)$ and $B=RC(Y,Y_0)$.
Let $x\in X$. Then
$$(\tcy0\circ f)(x)=\tcy0(f(x))=\s_{f(x)}^B=\{F\in B\st f(x)\in B\}.$$
Further,
$\fs(\tcx0(x))=\fs( \s_x^A)=\p_f\inv(\s_x^A)$. Set $\GA\ap= \p_f\inv(\s_x^A)$.
Then $\GA\ap=\{G\in B\st \p_f(G)\in\s_x^X\}=\{\cl_Y(G_0)\st G_0\in CO(Y_0), x\in\cl_X(X_0\cap f\inv(G_0))\}$.
Let $G_0\in CO(Y_0)$ and $\cl_Y(G_0)\in\GA\ap$. Then $x\in\cl_X(X_0\cap f\inv(G_0))$ and thus $f(x)\in f(\cl_X(X_0\cap f\inv(G_0))\sbe \cl_Y(f(X_0\cap f\inv(G_0)))\sbe \cl_Y(G_0)$. Therefore, $$\GA\ap\sbe\s_{f(x)}^B.$$ We have that $\GA\ap\in Grills(B)$. Hence $\GA\ap_r=r_{Y_0,Y}(\GA\ap)\in Grills(CO(Y_0))$. Thus, by (S2S4), there exists $y\in Y$ such that $\GA\ap_r=\GA_{y,Y_0}$. Then $$\GA\ap=\s_y^B.$$
Since $\GA\ap\sbe\s_{f(x)}^B$, we get, by Lemma \ref{pro41dv1}, that $\nu_{f(x)}^B\sbe\GA\ap$. According to \cite[Proposition 4.1]{DV1}, $\nu_{f(x)}^B$ is a filter of $B$. Hence, by Lemma \ref{grilllemma}, there exists an ultrafilter $u$ of $B$ such that $\nu_{f(x)}^B\sbe u\sbe \GA\ap$. Then $u\sbe \s_y^B$ and since $u$ is a grill of $B$, Lemma \ref{pro41dv1} implies that $\nu_y^B\sbe u$. So, we obtained that $\nu_{f(x)}^B\cup\nu_y^B\sbe u\sbe \GA\ap$. Then, for every $F\ap\in\nu_{f(x)}^B$ and every $G\ap\in\nu_y^B$, we have that $F\ap . G\ap\neq 0$, i.e. $\cl_Y(\int_Y(F\ap\cap G\ap))\nes$. Hence $\int_Y(F\ap\cap G\ap)\nes$ and thus $\int_Y(F\ap)\cap \int_Y(G\ap)\nes$, for every $F\ap\in\nu_{f(x)}^B$ and every $G\ap\in\nu_y^B$. Since $Y$ is a $T_0$-space, using Lemma \ref{lmconst}, we get that $y=f(x)$. Therefore
$\GA\ap=\s_{f(x)}^B$. Thus $\fs\circ \tcx0=\tcy0\circ f$ and hence $t$ is a natural isomorphism.

Finally, we will prove that
$$g: Id_{\Boo}\lra  D^t\circ  D^a,\mbox{ where } g(A)=g_A, \ \ \fa A\in\card\Boo,$$
%
 is a natural
isomorphism.

Let $\p\in\Boo(A,B)$ and $\ps=D^t(D^a(\p))$. We have
to prove that $g_{B}\circ\p=\ps\circ g_A$. Set $f=
D^a(\p)$, $(X,X_0)=D^a(A)$ and $(Y,Y_0)=D^a(B)$. Then $\ps=D^t(f)(=\p_f)$. Let
$a\in A$. Then $g_B(\p(a))=\{\GA\ap\in Y\st \p(a)\in\GA\ap\}.$
Further,
using (\ref{g0acl}), we get that
$$g_B(\p(a))=\cl_Y(s_B(\p(a)))\ \ \mbox{ and }\ \  g_A(a)=\cl_X(s_A(a)).$$
Thus $$\ps(g_A(a))=\cl_Y(Y_0\cap f\inv(s_A(a))).$$
Let $u\ap\in Y_0\cap f\inv(s_A(a))$. Then $u\ap\in Ult(B)$ and $f(u\ap)\in s_A(a)$. Hence $\p\inv(u\ap)\in s_A(a)=\{u\in Ult(A)\st a\in u\}$. Thus $a\in\p\inv(u\ap)$, i.e. $\p(a)\in u\ap$. Therefore $u\ap\in s_B(\p(a))$. So, $Y_0\cap f\inv(s_A(a))\sbe s_B(\p(a))$. Conversely, let $u\ap\in s_B(\p(a))$. Then $u\ap\in Y_0$ and $\p(a)\in u\ap$. Hence $a\in\p\inv(u\ap)=f(u\ap)$. Thus $f(u\ap)\in s_A(a)$. Therefore, $u\ap\in Y_0\cap f\inv(s_A(a))$. So, $Y_0\cap f\inv(s_A(a))\spe s_B(\p(a))$ and we get that $Y_0\cap f\inv(s_A(a))= s_B(\p(a))$. Hence $\ps(g_A(a))=\cl_Y(s_B(\p(a)))=g_B(\p(a))$.
 So,
$g$ is a natural isomorphism.

We have proved that $(D^t,D^a,g,t)$ is a duality between the categories $\ECS$ and $\Boo$.
 \sqs

\begin{defi}\label{ecscat}
\rm
(a) Let $\ECC$  be the category whose objects are all extremally connected spaces and whose morphisms are all continuous maps between the objects of $\ECC$ which preserve u-points (i.e., for every $X,Y\in|\ECC|$, $f\in\ECC(X,Y)$ iff $f$ is a continuous map and for every u-point $x\in X$, $f(x)$ is an u-point of $Y$).

\smallskip

\noindent(b) Let
$\CBool$ be the full subcategory of the category $\Boo$, whose objects are all complete Boolean algebras.
\end{defi}

\begin{rem}\label{remecscat}
\rm
(a) Clearly, $\ECC$ is indeed a category;

\smallskip

\noindent(b) Note that, according to Proposition \ref{proupointopen}, every open map between two objects of the category $\ECC$ is an $\ECC$-morphism.
\end{rem}

\begin{theorem}\label{thecscat}
The categories $\CBool$ and $\ECC$ are dually equivalent.
\end{theorem}

\doc Let $\CStone$ be the full subcategory of the category $\ECS$, whose objects are all Stone 2-spaces $(X,X_0)$ for which $X_0$ is extremally disconnected. We will first show that the categories $\CStone$ and $\ECC$ are isomorphic. Let us define two (covariant) functors
$$E_1:\ECC\lra\CStone\ \  \mbox{ and }\ \ E_2:\CStone\lra\ECC.$$
Let $X\in|\ECC|$ and $X_0=\{x\in X\st x$ is an u-point of $X\}$. Then, by Corollary \ref{grsp2grsp}, $(X,X_0)\in|\CStone|$ and we set $$E_1(X)=(X,X_0).$$
Let $(X,X_0)\in|\CStone|$. Then, by Lemma \ref{2grgr}, $X\in|\ECC|$ and we set $$E_2(X,X_0)=X.$$
Let $f\in\ECC(X,Y)$, $E_1(X)=(X,X_0)$ and $E_1(Y)=(Y,Y_0)$. Then, by the corresponding definitions, we get that $f$ is continuous and $f(X_0)\sbe Y_0$. Hence $f\in\CStone(E_1(X),E_1(Y))$ and we set $$E_1(f)=f.$$
Let $f\in \CStone((X,X_0),(Y,Y_0))$. Then, by Lemma \ref{2grgr} and Corollary \ref{grsp2grsp}, we get that $X_0=\{x\in X\st x$ is an u-point of $X\}$ and $Y_0=\{y\in Y\st y$ is an u-point of $Y\}$. Since $f(X_0)\sbe Y_0$, we obtain that $f\in\ECC(X,Y)=\ECC(E_2(X,X_0),E_2(Y,Y_0))$ and we set $$E_2(f)=f.$$
Obviously, $E_1$ and $E_2$ are functors. If $X\in|\ECC|$ then $E_2(E_1(X))=E_2(X,X_0)=X$. If $(X,X_0)\in|\CStone|$ then $E_1(E_2(X,X_0))=E_1(X)$. Using again Lemma \ref{2grgr} and Corollary \ref{grsp2grsp}, we get that $E_1(X)=(X,X_0)$. Hence $E_1(E_2(X,X_0))=(X,X_0)$. Now it becomes obvious that $E_1\circ E_2=Id_{\CStone}$ and $E_2\circ E_1=Id_{\ECC}$. So, the categories $\ECC$ and $\CStone$ are isomorphic. Let EDStone be the class of all extremally disconnected Stone spaces. Then, using the Stone Theorem that $S(|\CBool|)=$ EDStone, we get that the restrictions $D^a_{|\CBool}$ and $D^t_{|\CStone}$ of the duality functors $D^a$ and $D^t$ defined in the proof of Theorem \ref{constoneduality}, are duality functors between the categories $\CBool$ and $\CStone$.
Setting
\begin{equation}\label{dacdtc}
  D^a_c=E_2\circ D^a_{|\CBool}\ \ \mbox{ and }\ \ D^t_c=D^t_{|\CStone}\circ E_1,
\end{equation}
we obtain that
 $$ D^a_c:\CBool\lra\ECC\ \ \mbox{ and }\ \ D^t_c:\ECC\lra\CBool$$
are duality functors.
\sqs

\section{On a class of compact $T_0$ extensions}

\begin{defi}\label{niDE}
\rm
An {\em extension} of a space $X$ is a pair $(Y,f)$, where $Y$ is
a space and $f:X\lra Y$ is a dense embedding of $X$ into  $Y$.

Two extensions $(Y_i,f_i),\ i=1,2$, of $X$ are called {\em isomorphic}
(or {\em equivalent}\/)
if there exists a homeomorphism $\p:Y_1\lra Y_2$ such that $\p\circ f_1=f_2$.
Clearly, the relation of isomorphism   is an equivalence in  the class of
all extensions of $X$; the equivalence class of an extension $(Y,f)$ of $X$ will be denoted by $[(Y,f)]$.

We write $$(Y_1,f_1)\le (Y_2,f_2)$$
and say that the extension $(Y_2,f_2)$ is {\em projectively larger} than the
extension $(Y_1,f_1)$
if there exists a continuous mapping
$f:Y_2\lra Y_1$ such that $f\circ f_2=f_1$. This relation is a
{\em preorder}\/ (i.e., it is reflexive and transitive).
Setting for every two  extensions $(Y_i,f_i),\ i=1,2$, of a  space $X$,
$[(Y_1,f_1)]\le [(Y_2,f_2)]$  iff $(Y_1,e_1)\le (Y_2,e_2),$ we obtain a well-defined relation on the class of all, up to equivalence,
extensions of $X$; obviously, it is also a preorder (see, e.g., \cite{Ba}).

We write $$(Y_1,f_1)\le_{in} (Y_2,f_2)$$
and say that the extension $(Y_2,f_2)$ is {\em injectively larger} than the
extension $(Y_1,f_1)$ if there exists a continuous mapping
$f:Y_1\lra Y_2$ such that $f\circ f_1=f_2$
and $f$ is a homeomorphism from $Y_1$ to the subspace $f(Y_1)$
of $Y_2$. This relation is  a preorder.
Setting for every two  extensions $(Y_i,f_i),\ i=1,2$, of a  space $X$,
$[(Y_1,f_1)]\le_{in} [(Y_2,f_2)]$  iff $(Y_1,e_1)\le_{in} (Y_2,e_2),$ we obtain a well-defined relation on the class of all, up to equivalence,
extensions of $X$; obviously, it is also a preorder (see, e.g., \cite{Ba}).
\end{defi}

\begin{nota}\label{csrextord}
\rm
Let $Y$ be a space. We will denote by $\CSR(Y)$ (resp., by $\CCSR(Y)$) the class of all, up  to equivalence, (connected) C-semiregular extensions of $Y$.

Recall that if $B$ is a Boolean algebra, then we  denote by $\CRel(B)$ (resp., $\CCRel$) the  set of all (connected) contact relations on  $B$. We define a relation $``\le$" on the set $\CRel(B)$ setting, for any $C_1,C_2\in\CRel(B)$,
$C_1\le C_2\iff C_1\spe C_2$. We will denote again by $``\le$" the restriction of the relation $``\le$" to the set $\CCRel$.
\end{nota}

\begin{theorem}\label{csemregext}
Let $Y$ be an extremally disconnected compact Hausdorff space and $B=RC(Y)$. Then the ordered sets $(\CRel(B),\le)$ and $(\CSR(Y),\le)$, as well as the ordered sets $(\CRel(B),\sbe)$ and $(\CSR(Y),\le_{in})$, are isomorphic (see Definition \ref{niDE} for the relations $``\le$" and $``\le_{in}$" on $\CSR(Y)$). Also, the ordered sets
$(\CCRel(B),\le)$ and $(\CCSR(Y),\le)$, as well as the ordered sets $(\CCRel(B),\sbe)$ and $(\CCSR(Y),\le_{in})$, are isomorphic.
\end{theorem}

\doc Let $(X,f)$ be a C-semiregular extensions of $Y$.
Set $X\ap=f(Y)$. Then,
 clearly,
 the map $$e:(RC(X\ap),\d_{(X,X\ap)})\lra (RC(X),C_X), \ \ F\mapsto\cl_X(F),$$
is a CA-isomorphism (note that $RC(X\ap)=CO(X\ap)$). For every $F,G\in B$, set
\begin{equation}\label{fcxg}
FC_{(X,f)}G\iff \cl_X(f(F))\cap\cl_X(f(G))\nes,
\end{equation}
 i.e., $FC_{(X,f)}G\iff f(F)\d_{(X,X\ap)} f(G)$. Then, obviously, $(B,C_{(X,f)})$ is a contact algebra.
 Set $$\p(X,f)=(B,C_{(X,f)}).$$
 Clearly, two equivalent extension of $Y$ define two coinciding contact relations on the Boolean algebra $B$.
 Thus we have that  $\p([(X,f)])=(B,C_{(X,f)})$ and, for simplicity, we will denote by the same letter $\p$ the induced  map on the set of equivalence classes of the C-semiregular extensions of $Y$.


Conversely, let $C$ be a contact relation on the Boolean algebra $B$ and let $(\wh{X},\wh{X}_0)$ be the canonical 2-contact space of the complete contact algebra $(B,C)$ (see Definition \ref{contactspace}(c)). Then, by the definition of the space $\wh{X}_0$ and the Stone Representation Theorem, we have that the map
$$\wh{f}:Y\lra \wh{X},\ \  y\mapsto u_y,$$
(see (\ref{ux}) for the notation $u_y$) is a homeomorphic embedding and $\wh{f}(Y)=\wh{X}_0$. Hence, $(\wh{X},\wh{f})$ is an extension of the space $Y$. Using Lemma \ref{csem2cont}, we get that $\wh{X}$ is a C-semiregular space. So, $(\wh{X},\wh{f})$ is a C-semiregular extension of the space $Y$. Set $$\psi(B,C)=(\wh{X},\wh{f}).$$

Let $(X,f_0)$ be a C-semiregular extensions of $Y$,  $(B,C)=\p(X,f_0)$ and $(\wh{X},\wh{f}_0)=\psi(B,C)$. We will show that $(X,f_0)$ and $(\wh{X},\wh{f}_0)$ are isomorphic extensions of $Y$. As we have already seen, the map
$$e:(RC(f_0(Y)),\d_{(X,f_0(Y))})\lra(RC(X),C_X),\ \ G\mapsto \cl_X(G),\ \  \mbox{ is a CA-isomorphism.}$$
Clearly, the map
$$\g_{f_0}^0:(B,C)\lra (RC(f_0(Y)),\d_{(X,f_0(Y))}), \ \  F\mapsto f_0(F),\ \  \mbox{ is a CA-isomorphism.}$$
Set $\g^0=e\circ\g_{f_0}^0$. Then
$$\g^0:(B,C)\lra (RC(X),C_X),\ \  F\mapsto\cl_X(f_0(F)),\ \  \mbox{ is a CA-isomorphism.}$$
Thus the map
$$\g\ap:Clans(B,C)\lra Clans(RC(X),C_X),\ \  \GA\mapsto\g^0(\GA),\ \  \mbox{ is a bijection}.$$
Since $X$ is a C-semiregular space, Lemma \ref{dv141ii42i} implies that the map
$$\kappa: X\lra Clans(RC(X),C_X),\ \  x\mapsto\s_x,\ \ \mbox{ is a bijection.}$$
Let $\l=\k\inv$. Then we get that the map
$$f:\wh{X}\lra X,\ \  \GA\mapsto\l(\g\ap(\GA)),\ \ \mbox{ is a bijection}$$
(i.e., we set $f=\l\circ\g\ap$). We will show that $f\circ\wh{f}_0=f_0$. Indeed, let $y\in Y$. Then $\wh{f_0}(y)=u_y\in\wh{X}$ and $f(\wh{f_0}(y))=f(u_y)=\l(\g\ap(u_y))$. Set $x=\l(\g\ap(u_y))$. Then
$\k(x)=\g^0(u_y)$, i.e., $\s_x=\{\cl_X(f_0(F))\st F\in B, y\in F\}=\{\cl_X(f_0(F))\st F\in B, f_0(y)\in f_0(F)\}=\{\cl_X(G)\st G\in RC(f_0(Y)), f_0(y)\in G\}=e(u_{f_0(y)})$. Hence
\begin{equation}\label{euxx0}
e(u_{f_0(y)})=\s_x\ \ \mbox{ and, thus, }\ \ r(\s_x)=u_{f_0(y)}.
 \end{equation}
  Suppose that $x\neq f_0(y)$. Since $X$ is $T_0$ and semiregular, we get that there exists $F\in RC(X)$ such that $|F\cap\{x,f_0(y)\}|=1$. If $x\in F$, then $f_0(y)\nin F$. Thus $F\in\s_x$ and $F\cap f_0(Y)\nin u_{f_0(y)}$. Since $r(F)=F\cap f_0(Y)$, we get a contradiction (see (\ref{euxx0})). If $f_0(y)\in F$, then $x\nin F$. Thus $F\nin\s_x$ and $F\cap f_0(Y)\in u_{f_0(y)}$. Since $e(F\cap f_0(Y))=F$ (by Lemma \ref{isombool}), we get  a contradiction (see again (\ref{euxx0})). Hence $x=f_0(y)$.  Therefore,
$$f\circ\wh{f}_0=f_0.$$
We will now show that $f$ is a homeomorphism. Let $F\in B$. Then, using the notation of Definition \ref{canspa}, we obtain that $f(g_B(F))=f(\{\GA\in\wh{X}\st F\in\GA\})=\{f(\GA)\st F\in\GA\}=\{f(\GA)\st \g^0(F)\in\g^0(\GA)\}=\{f(\GA)\st \cl_X(f_0(F))\in\g\ap(\GA)\}=\{\l(\g\ap(\GA))\st \cl_X(f_0(F))\in\k(\l(\g\ap(\GA)))\}=\{f(\GA)\st \cl_X(f_0(F))\in\s_{f(\GA)}\}=\{f(\GA)\st f(\GA)\in\cl_X(f_0(F))=\cl_X(f_0(F))$. So, $f(g_B(F))=\cl_X(f_0(F))$, for every $F\in B$.  Since $f$ is a bijection, we also get that
for every $F\in B$, $f\inv(\cl_X(f_0(F)))=g_B(F)$. Now, using the fact that $\{\cl_X(f_0(F))\st F\in B\}=RC(X)$ and that $RC(X)$ and $\{g_B(F)\ F\in B\}$ are closed bases of, respectively, $X$ and $\wh{X}$, we get that $f$ is a homeomorphism. Therefore,
$$\psi(\p((X,f_0)))\  \mbox{ is isomorphic to }\ (X,f_0).$$

Let now $C$ be a contact relation on the Boolean algebra $B$, $\psi(B,C)=(\wh{X},\wh{f})$ and $\p(\psi(B,C))=(B,\wh{C})$. We will show that $C\equiv \wh{C}$. We have that for every $F,G\in B$, $F\wh{C}G\iff \cl_{\wh{X}}(\wh{f}(F))\cap\cl_{\wh{X}}(\wh{f}(G))\nes$. Recall that the set $\{h_B(H)\st H\in B\}$, where $h_B(H)=\{\GA\in \wh{X}\st H\nin\GA\}$, is an open base of $\wh{X}$. Let us show that if $H\in B$ and $\GA\in\wh{X}$, then
\begin{equation}\label{gaclfh}
\GA\in\cl_{\wh{X}}(\wh{f}(H))\iff H\in\GA.
 \end{equation}
 Indeed, using the fact that $\GA$ satisfies condition (Clan2) (see Definition \ref{clandef}), we get that  ($\GA\in\cl_{\wh{X}}(\wh{f}(H))$) $\iff$ (for every $P\in B\stm\GA$, $h_B(P)\cap\{u_y\st y\in H\}\nes$) $\iff$ (for every $P\in B\stm\GA$, there exists $y\in H$ such that $P\nin u_y$) $\iff$ (for every $P\in B\stm\GA$, there exists $y\in H$ such that $y\nin P$) $\iff$ (for every $P\in B\stm\GA$, $H\not\sbe P$) $\iff$ ($H\in\GA$).
 So, (\ref{gaclfh}) is verified. Now, we get, using (\ref{gaclfh}) and Lemma \ref{lmpc}(c), that for every $F,G\in B$, $F\wh{C}G\iff \ex \GA\in(\cl_{\wh{X}}(\wh{f}(F))\cap\cl_{\wh{X}}(\wh{f}(G)))\iff (\ex\GA\in\wh{X})(F,G\in\GA)\iff FCG$. Therefore $C\equiv \wh{C}$.

 So, the correspondence $\p$  is a bijection between the  set of all, up  to equivalence, C-semiregular extensions of $Y$  and the  set of all contact relations on the Boolean algebra $B$. Let us show that $\p$ is an isomorphism.

 Let $(X_i,f_i)$, $i=1,2$, be two C-semiregular extensions  of $Y$ and $[(X_1,f_1)]\ge[(X_2,f_2)]$ or $[(X_1,f_1)]\le_{in}[(X_2,f_2)]$. Then there exists a continuous mapping $f:X_1\lra X_2$ such that $f\circ f_1=f_2$. Set $(B,C_i)=\p(X_i,f_i)$, $i=1,2$. Let $F,G\in B$ and $F(-C_2)G$. Then, by (\ref{fcxg}), $\cl_{X_2}(f_2(F))\cap\cl_{X_2}(f_2(G))=\ems$.
 Hence $f\inv(\cl_{X_2}(f_2(F)))\cap f\inv(\cl_{X_2}(f_2(G)))=\ems$. Since $f(\cl_{X_1}(f_1(F)))\sbe \cl_{X_2}(f(f_1(F)))=\cl_{X_2}(f_2(F))$, we get that
 $\cl_{X_1}(f_1(F))\sbe f\inv(\cl_{X_2}(f_2(F)))$. Analogously, $\cl_{X_1}(f_1(G))\sbe f\inv(\cl_{X_2}(f_2(G)))$. Thus $\cl_{X_1}(f_1(F))\cap\cl_{X_1}(f_1(G))=\ems$. Using once more (\ref{fcxg}), we get that $F(-C_1)G$. Therefore $C_1\sbe C_2$, i.e., $C_1\ge C_2$.

 Conversely, let $C_1$ and $C_2$ be two contact relations on $B$ and $C_1\ge C_2$, i.e., $C_1\sbe C_2$. Set $(X_i,f_i)=\psi(B,C_i)$, $i=1,2$. Then $\p\inv(B,C_i)=[(X_i,f_i)]$. By the definition of the map $\psi$, we have that for $i=1,2$, $X_i=Clans(B,C_i)$, the topology on $X_i$ is generated by the closed base $\{\{\GA\in Clans(B,C_i)\st F\in \GA\}\st F\in B\}$, and $f_i(y)=u_y$, for every $y\in Y$. Since $C_1\sbe C_2$, we get that $Clans(B,C_1)\sbe Clans(B,C_2)$. Now we define
 $$f:X_1\lra X_2,\ \ \GA\mapsto\GA.$$ Then, for every $y\in Y$, $f(f_1(y))=f(u_y)=u_y=f_2(y)$. Hence, $f\circ f_1=f_2$. Since, for every $F\in B$, $f\inv(\{\GA\in X_2\st F\in\GA\})=\{\GA\in X_1\st F\in\GA\}$, we get that $f$ is a continuous map. Therefore, $[(X_1,f_1)]\ge[(X_2,f_2)].$ Let us note that $f$ is even an embedding. Indeed, $f$ is injective and for every $F\in B$, $f(\{\GA\in X_1\st F\in\GA\})=f(X_1)\cap\{\GA\in X_2\st F\in\GA\}$. Hence, if $f\ap:X_1\lra f(X_1)$ is the restriction of $f$ and $g\ap=(f\ap)\inv:f(X_1)\lra X_1$, then, for every $F\in B$, $(g\ap)\inv(\{\GA\in X_1\st F\in\GA\})=f(X_1)\cap\{\GA\in X_2\st F\in\GA\}$. So that, $f$ is an embedding. Hence, $[(X_1,f_1)]\le_{in}[(X_2,f_2)].$

 Therefore, $\p$ is an isomorphism between the ordered sets  $(\CSR(Y),\le)$ and $(\CRel(B),\le)$, and also between the ordered sets $(\CSR(Y),\le_{in})$ and $(\CRel(B),\sbe)$. Clearly, this implies that $\CSR(Y)$ is a set and the  preorders $``\le$" and $``\le_{in}$" on $\CSR(Y)$, defined in Definition \ref{niDE}, are, in fact,  orders.

 Now, the assertions about connected contact relations on $B$ follow immediately.
 \sqs

  Noting that the Stone space of a Boolean algebra $B$ is extremally disconnected iff $B$ is complete (see, e.g., \cite{kop89}), the above theorem can be reformulated as follows:

 \begin{theorem}\label{contrelcomplba}
 Let $B$ be a complete Boolean algebra and $Y=S(B)$ be its Stone space. Then the ordered sets $(\CRel(B),\le)$ and $(\CSR(Y),\le)$, as well as the ordered sets $(\CRel(B),\sbe)$ and $(\CSR(Y),\le_{in})$, are isomorphic (see Definition \ref{niDE} for the relations $``\le$" and $``\le_{in}$" on $\CSR(Y)$). Also, the ordered sets
$(\CCRel(B),\le)$ and $(\CCSR(Y),\le)$, as well as the ordered sets $(\CCRel(B),\sbe)$ and $(\CCSR(Y),\le_{in})$, are isomorphic.
 \end{theorem}

 \begin{cor}\label{grsmcsrgext}
 Let $Y$ be an extremally disconnected compact Hausdorff space and $B=RC(Y)$. Then the ordered sets  $(\CSR(Y),\le)$ and $(\CSR(Y),\le_{in})$ have  largest and  smallest elements. The largest (resp., the smallest) element of the ordered set $(\CSR(Y),\le)$ coincides with the smallest (resp., the largest) element of
 $(\CSR(Y),\le_{in})$. For the largest element $[(\g Y,\g_Y)]$ of the ordered set
 $(\CSR(Y),\le_{in})$, we have that $\g Y$ is an extremally connected space (in fact, $\g Y=D^a_c(B)$ (see (\ref{dacdtc}) for $D^a_c$)). Also, if $(cY,c)$ is a C-semiregular extension of $Y$ and $cY$ is an extremally connected space then the C-semiregular extensions $(cY,c)$ and $(\g Y,\g_Y)$ of $Y$ are equivalent.
 \end{cor}

 \doc Clearly, by Theorem \ref{csemregext}, the smallest (resp., the largest) element of the ordered set $(\CSR(Y),\le)$ is the largest (resp., the smallest) element of the ordered set $(\CSR(Y),\le_{in})$. So that we will regard only the ordered set $(\CSR(Y),\le_{in})$. By Example \ref{extrcr}, the ordered set $(\CRel(B),\sbe)$, where $B=RC(Y)(=CO(Y))$, has  largest and  smallest elements.  Thus, by Theorem \ref{csemregext}, the ordered set  $(\CSR(Y),\le_{in})$ also has  largest and  smallest elements. The fact that it has  smallest element follows also from Corollary \ref{novocor}: this is the equivalence class of the extension $(Y,id_Y)$ of $Y$. It is also obvious that it corresponds to the contact relation $\rho_s$ on $B$ (see the formula (\ref{fcxg})). For the largest element $[(\g Y,\g)]$ of the ordered set
 $(\CSR(Y),\le_{in})$, we have that the map $\g$ is defined by the formula $\g_Y(y)=u_y$, for every $y\in Y$ (see (\ref{ux}) for the notation $u_y$), and
 $\g Y=Clans(B,\rho_l)=Grills(B)=D^a_c(B)$. Thus $\g Y$ is an extremally connected space   (see Theorem \ref{thecscat}).

 Let $(cY,c)$ be a C-semiregular extension of $Y$ and let $cY$ be an extremally connected space.
 By Proposition \ref{grspCsemrs}(a), we have that the standard contact relation $C_{cY}$ on $RC(cY)$ coincides with the largest contact relation $\rho_l^{RC(cY)}$ on the Boolean algebra $RC(cY)$. Then, using (\ref{fcxg}), we obtain that the contact relation $C_{(cY,c)}$ on the Boolean algebra $B$, corresponding to the C-semiregular extension $(cY,c)$ of $Y$ (see the proof of Theorem \ref{csemregext}), coincides with the  largest contact relation $\rho_l^B$ on the Boolean algebra $B$. Therefore, $(cY,c)$ corresponds to $\rho_l^B$; thus $(cY,c)$ and $(\g Y,\g_Y)$  are equivalent C-semiregular extensions of $Y$.
 \sqs

\begin{theorem}\label{extcontextr}
 Let $X$ and $Y$ be two extremally disconnected compact Hausdorff spaces, $(cX,c)$ be an arbitrary C-semiregular extension of $X$ and $f:X\lra Y$ be a  continuous map.  Then there exists a continuous map $f\ap:cX\lra \g Y$ such that $\g_Y\circ f=f\ap\circ c$ (see  Corollary \ref{grsmcsrgext} for $(\g Y,\g_Y)$)  (i.e., supposing that $c$ and $\g_Y$ are the embedding maps of $X$ and $Y$ in, respectively, $cX$ and $\g Y$, we get that
 $f$  can be extended to a continuous map $f\ap:cX\lra \g Y$). In particular, every continuous map $f:X\lra Y$ can be $``$extended" to a continuous map $\g f:\g X\lra\g Y$ (i.e.,  $\g_Y\circ f=\g f\circ \g_X$).
\end{theorem}

\doc Since $(cX,c)\le_{in}(\g X,\g_X)$ (see Corollary \ref{grsmcsrgext}), we can regard $cX$ as a subspace of $\g X$. Thus, it is enough to prove only that there exists a continuous  map $\g f:\g X\lra\g Y$ such that  $\g_Y\circ f=\g f\circ \g_X$. Regard the Boolean algebras $A=RC(X)$ and $B=RC(Y)$. By the Stone Duality, the map $$\p_f=S(f):B\lra A, \ \  G\mapsto f\inv(G),$$ is a Boolean homomorphism. Hence, using Corollary \ref{grsmcsrgext} and Theorem \ref{thecscat}, we get that $D^a_c(\p_f):\g X\lra \g Y$ is a continuous map. Set $\g f=D^a_c(\p_f)$. We will show that $\g_Y\circ f=\g f\circ \g_X$. Let $x\in X$. Set $y=f(x)$. Using  Corollary \ref{grsmcsrgext},
(\ref{deftetaphic}) and (\ref{dacdtc}), we get that $\g f(\g_X(x))=\g f(u_x^A)=(\p_f)\inv(u_x^A)$ and $\g_Y(f(x))=\g_Y(y)=u_y^B$. So, we have to show that $u_y^B=(\p_f)\inv(u_x^A)$. Let $G\in B$. Then we have that $G\in (\p_f)\inv(u_x^A)\iff \p_f(G)\in u_x^A\iff f\inv(G)\in u_x^A\iff x\in f\inv(G)\iff f(x)\in G\iff y\in G$. Therefore, $u_y^B=(\p_f)\inv(u_x^A)$. Thus, $\g_Y\circ f=\g f\circ \g_X$.
\sqs

\begin{theorem}\label{extcontextrconn}
Let $X$ be an extremally disconnected compact Hausdorff space,  $Z$ be an extremally connected  space, $(cX,c)$ be an arbitrary C-semiregular extension of $X$ and $f:X\lra Z$ be a  continuous map such that, for every $x\in X$, $f(x)$ is an u-point of $Z$.  Then there exists a continuous map $f\ap:cX\lra Z$ such that $f=f\ap\circ c$   (i.e., supposing that $c$ is the embedding map of $X$ in $cX$, we get that
 $f$  can be extended to a continuous map $f\ap:cX\lra Z$). In particular, every open map $f:X\lra Z$ can be $``$extended" to a continuous map $f\ap:cX\lra Z$ (i.e., $f=f\ap\circ c$).
\end{theorem}

\doc Set $Y=\{z\in Z\st z$ is an u-point of $Z\}$. Then, by Corollary \ref{grsp2grsp}, $Y$ is a dense  extremally disconnected compact Hausdorff subspace of $Z$. Setting $i_Y:Y\lra Z$ to be the embedding of $Y$ in $Z$, we get (by Proposition \ref{grspCsemrs}(b)) that $(Z,i_Y)$ is a C-semiregular extension of $Y$. Moreover, Corollary \ref{grsmcsrgext} implies that $(Z,i_Y)$ and $(\g Y,\g_Y)$ are equivalent C-semiregular extensions of $Y$. Since $f(X)\sbe Y$, our assertion follows now from Theorem \ref{extcontextr}.

Finally, note that if $f$ is an open map, then, by    Proposition \ref{remupoint}(a)  and Proposition \ref{proupointopen}, we have that, for every $x\in X$, $f(x)$ is an u-point of $Z$ .
\sqs

\end{document}